\documentclass[reqno]{amsart}
\usepackage{amsmath, amsthm, amssymb, amstext}

\usepackage{orcidlink}
\usepackage{xcolor}
\definecolor{revisiongreen}{RGB}{0,115,35}
\newcommand{\newadd}[1]{{\color{revisiongreen}#1}}
\usepackage{hyperref}
\newcommand{\rev}[1]{{\color{revisiongreen}#1}}
\newcommand{\Pohozaev}{\rev{\mbox{Poho\v{z}aev}}}
\hypersetup{
	pdfborder={0 0 0},
	colorlinks=false,
	linkcolor=black,
	urlcolor=black,
	citecolor=black
}

\usepackage{enumitem}
\allowdisplaybreaks
\usepackage[toc,page]{appendix}
\usepackage{graphicx}
\usepackage{caption}
\usepackage{mathtools}
\usepackage{tikz}

\newtheorem{theorem}{Theorem}[section]
\newtheorem{remark}[theorem]{Remark}
\newtheorem{proposition}[theorem]{Proposition}
\newtheorem{corollary}[theorem]{Corollary}
\newtheorem{definition}[theorem]{Definition}

\newtheorem{lemma}[theorem]{Lemma}

\newcommand{\R}{\mathbb{R}}

\numberwithin{theorem}{section}
\numberwithin{equation}{section}
\renewcommand{\theequation}{\arabic{section}.\arabic{equation}}
\renewcommand{\thetheorem}{\arabic{section}.\arabic{theorem}}
\title[Fractional Logarithmic Laplacian]{Normalized Solutions to Nonlinear Schrödinger Equations Involving the Fractional-Logarithmic Laplacian}

\date{\today}

\author[A. Gouasmia]{Abdelhamid Gouasmia}

\address{\textbf{Abdelhamid Gouasmia}\\
Department of Mathematics and Computer Science\\
Faculty of Exact Sciences and Natural and Life Sciences\\
Echahid Cheikh Larbi Tebessi University, T\'ebessa 12000, Algeria\\
Laboratoire d'\'Equations aux D\'eriv\'ees Partielles Non Lin\'eaires et Histoire des Math\'ematiques, ENS, B.P. 92, Vieux Kouba, 16050 Algiers, Algeria}
\email{gouasmia.abdelhamid@gmail.com}

\author[H. Hajaiej]{Hichem Hajaiej}

\address{
		\textbf{Hichem Hajaiej}\\
		Department of Mathematics, California State University, Los Angeles, CA 90032, USA.}
\email{hichem.hajaiej@gmail.com}


\usepackage{mathrsfs}
\newcommand{\E}{\mathcal E_p}
\newcommand{\Q}{\mathcal Q}
\newcommand{\HH}{H^{s+\log}(\R^N)}
\newcommand{\Sc}{\mathcal S_c}
\newcommand{\Pc}{\mathcal P_c}
\newcommand{\Bc}{\mathbf S_c}
\newcommand{\Js}{\mathcal J_s}
\newcommand{\Jp}{\mathcal J_+}
\newcommand{\Jm}{\mathcal J_-}
\newtheorem*{fibercriterion}{Fiber-root criterion}

\begin{document}

\begin{abstract}
We study real normalized standing waves for a focusing nonlinear Schr\"odinger equation with fractional--logarithmic dispersion. The operator has Fourier symbol $|\xi|^{2s}\ln|\xi|^2$ and spectral bottom $-1/(es)$. In the fractional mass-subcritical and mass-critical ranges, $2<p\leq2+4s/N$, the physical energy
\[
 I(u)=\tfrac12\langle(-\Delta)^{s+\log}u,u\rangle-\tfrac1p\|u\|_p^p
\]
attains its minimum at every prescribed $L^2$-norm $c>0$, with
\[
 \inf_{\|u\|_2=c}I(u)<-\frac{c^2}{2es}.
\]
This strict comparison excludes vanishing and yields strict subadditivity and compactness modulo translations. At the fractional critical power, the logarithmic strengthening of the kinetic energy removes any small-mass restriction. In the mass-supercritical, Sobolev-subcritical range, we obtain positive-energy minimizers on the upper Poho\u{z}aev branch at sufficiently small prescribed norm. For $N\geq2$ and $2+4s/N<p\leq2_s^*$, a separate radial minimax construction yields infinitely many normalized solutions at every prescribed norm, including the Sobolev-critical endpoint, with energies and multipliers tending to infinity. For the subcritical and critical minimizers, we establish small-mass energy and multiplier limits, quantitative gap bounds with matching powers in dimension one, concentration near the nonzero Fourier sphere $|\xi|=e^{-1/(2s)}$, and spatial delocalization of the mass-normalized profiles. Anisotropic oscillatory packets strengthen the lower gap bounds in higher dimensions, and the normalized profiles vanish in every admissible $L^q$, including the finite Sobolev endpoint.
\end{abstract}

\setcounter{tocdepth}{1}
\maketitle
\tableofcontents

{

\section{Introduction}
\label{section1}

The study of nonlocal dispersive equations has attracted considerable
attention in recent years, particularly in connection with the
fractional Laplacian and its applications in various areas of
mathematics and physics. For $s\in(0,1)$, the fractional Laplace operator
is defined by

\begin{align}\label{equ123}
(-\Delta)^s u(x)
:= C(N,s)\,\mathrm{P.V.}\!\int_{\mathbb{R}^N}
\frac{u(x)-u(y)}{|x-y|^{N+2s}}\,dy,
\end{align}
where $\mathrm{P.V.}$ denotes the Cauchy principal value and
$C(N,s)$ is the normalization constant given by
\begin{equation}\label{equ00}
C(N,s)
:=
\frac{s\,2^{2s}\Gamma\!\left((N+2s)/2\right)}
{\pi^{\frac{N}{2}}\Gamma(1-s)},
\end{equation}
with $\Gamma$ denoting the Gamma function. For further background
on fractional operators and their fundamental properties, we refer
the reader to
\cite{Bertoin1996,DiNezza-Palatucci-Valdinoci2012,Vazquez2012,Laskin2000}.

In particular, the fractional Schrödinger equation has received
considerable attention in the physics literature, owing to its
natural emergence in a variety of scientific and engineering
applications.  Its dispersive operator is closely related to the generators of jump-type L\'evy processes in probability theory. In this context, the models
introduced in \cite{Laskin2000,Laskin2002} are particularly
noteworthy, as they formulate the Feynman path integral in terms
of L\'evy flight trajectories rather than classical Brownian motion
paths. Applications of fractional Schr\"odinger dispersion to optics are discussed in \cite{Longhi2015}; see also \cite{GuoXu2006} for further physical models. From a mathematical perspective, fractional Schrödinger equations
provide a framework for modeling nonlocal dispersive phenomena.
Motivated by their broad range of applications and their intrinsic
analytical interest, these equations have been the subject of
substantial research. In particular, we refer to
\cite{Cho2014, Cho2013, CarlesDong2024} for results on the Cauchy
problem, to \cite{LuoHajaiej2022, CherguiGouHajaiej2023} for
the existence of normalized solutions, and to
\cite{Hajaiej2013, Hajaiej2014} for studies of the symmetry of
minimizers. This list is by no means exhaustive.

Recent work on fractional models with small order $s$ has highlighted applications of the operator in \eqref{equ123}; see \cite{Antil-Bartels,Pellacci-Verzini,Feulefack-Jarohs}. These advances have prompted the introduction of new classes of nonlocal operators, particularly those associated with logarithmic-type Laplacians (see \cite[Theorem 1.1]{Chen-Weth} for the linear setting and \cite[Theorem 1.1]{Dyda-Jarohs-Sk} for the nonlinear framework). The growing interest in these operators is largely motivated by their wide range of applications. 
In this regard, \cite{Chen-Weth,Dyda-Jarohs-Sk} develop Dirichlet theories for the logarithmic Laplacian and logarithmic $p$-Laplacian, respectively. These order-zero operators must be distinguished from the fractional--logarithmic operator at a fixed $s>0$.
 More recently, considerable attention has been devoted to
differentiating the fractional Laplace operator, not only in
the limit as $s\to0^+$, but also with respect to the fractional
parameter at an arbitrary order $s\in(0,1)$. In this direction,
\cite{Chen-Chen-Hauer} introduced the fractional logarithmic
Laplacian in the linear setting. More precisely, for $u\in C_c^2(\mathbb R^N)$ and $x\in\mathbb R^N$, the fractional logarithmic Laplacian $(-\Delta)^{s+\log}$ is defined by

\begin{equation*}
\begin{aligned}
&(-\Delta)^{s+\log}u(x)
:=\left.\frac{d}{dt}(-\Delta)^t u(x)\right|_{t=s}\\[4pt]
&\qquad =B(N,s)(-\Delta)^s u(x) +C(N,s)\,\mathrm{P.V.}\!\int_{\mathbb{R}^N}
\frac{u(x)-u(y)}{|x-y|^{N+2s}}
\bigl(-2\ln|x-y|\bigr)\,dy,
\end{aligned}
\end{equation*}
where $C(N,s)$ is defined in \eqref{equ00}, and
\begin{equation*}
B(N,s)
:=\frac{d}{ds}\ln C(N,s)
=2\ln 2+\frac{1}{s}
+\psi(1-s)
+\psi\!\left(\frac{N+2s}{2}\right),
\end{equation*}
with $\psi$ denoting the digamma function. For the corresponding
nonlinear setting, namely the fractional logarithmic
$q$-Laplacian $(-\Delta)_q^{s+\log}$, with $q>1$, we refer
the reader to \cite{BahrouniGouasmiaHajaiejOuannasser2026}.
In these works, the associated functional framework is
developed by introducing suitable energy spaces and establishing
embedding results. Within this framework, several problems
involving these operators have been investigated.

We investigate normalized standing waves for fractional--logarithmic dispersion. The logarithmic nonlinearity studied in \cite{CarlesDong2024} and the order-zero heat equation studied in \cite{ChenVeron2024} concern different equations. We make no claim here about well-posedness of the evolution problem. The stationary problem is formally associated with the following focusing evolution equation:

\begin{equation}\tag{P}\label{equ2}
i\,\partial_t\Psi(x,t)-(-\Delta)^{s+\log}\Psi(x,t)
=-|\Psi(x,t)|^{p-2}\Psi(x,t),
\quad (x,t)\in\mathbb{R}^N\times(0,\infty),
\end{equation}
subject to the initial condition
\begin{equation}\tag{I}\label{equ22}
\Psi(x,0)=\Psi_0(x),
\qquad x\in\mathbb{R}^N.
\end{equation}
Here, $p\in(2,\infty)$ is the exponent characterizing the nonlinear interaction, while the logarithmic character of the equation is captured by the operator $(-\Delta)^{s+\log}$.

The mass has different physical interpretations depending on the model.
In nonlinear optics, it corresponds to the total beam power;
in quantum mechanics, it is the integral of the probability density;
and in Bose--Einstein condensation, it represents the total
number of particles. These interpretations naturally motivate
the investigation of solutions preserving the qualitative
features of the evolution, particularly standing wave solutions. To this end, we seek solutions of \eqref{equ2} of the form
\begin{equation}\label{equ1}
\Psi(x,t)=e^{i\lambda t}u(x),
\qquad \lambda\in\mathbb{R},
\end{equation}
where $\lambda$ denotes the frequency parameter. For sufficiently
regular solutions to \eqref{equ2}--\eqref{equ22}, the mass is
conserved along the flow, that is,
\[
\|\Psi(t)\|_{L^2(\mathbb{R}^N)}
=
\|\Psi_0\|_{L^2(\mathbb{R}^N)},
\quad \text{for all }t\in[0,T).
\]
Substituting the ansatz \eqref{equ1} into \eqref{equ2} yields the stationary problem:
\begin{equation}\label{equ3}
(-\Delta)^{s+\log}u+\lambda u
=
|u|^{p-2}u
\quad\text{in }\mathbb{R}^N.
\end{equation}
We study standing waves satisfying a prescribed $L^2$-norm constraint. Such solutions are
commonly referred to as normalized solutions. More precisely,
for a given $c>0$, we seek solutions satisfying

\begin{equation}\label{equ5}
\int_{\mathbb{R}^N}|u|^2\,dx=c^2.
\end{equation}

Our objective is to establish the existence of normalized solutions of \eqref{equ3}. Here $c$ denotes the prescribed $L^2$-norm, so the mass is $c^2$. In this setting, normalized solutions are
characterized as critical points of the associated energy
functional restricted to the $L^2$-sphere. The parameter $\lambda\in\mathbb R$ arises as a Lagrange multiplier and is not prescribed a priori.
This variational formulation is particularly meaningful
from a physical viewpoint, since the $L^2$-norm is conserved
along the evolution flow. Moreover, such a characterization
is essential in the investigation of the orbital stability
or instability of solutions; see
\cite{Cho2014,LuoHajaiej2022,CherguiGouHajaiej2023,Soave2020,
	JeanjeanLuo2013} and the references therein. A major difficulty in the present study stems from the specific structure of the fractional logarithmic operator,
which substantially distinguishes the problem from the
standard fractional framework. Indeed, the kernel associated
with $(-\Delta)^{s+\log}$ is given by
\[
\mathbf{K}^{s+\log}(r)
:=
\frac{B(N,s)C(N,s)-2C(N,s)\ln r}
{r^{N+2s}},
\qquad r>0.
\]
Unlike the standard fractional kernel, this logarithmic
kernel changes sign at large distances and exhibits an
additional logarithmic singularity near the diagonal.
Consequently, the associated energy functional is not
necessarily positive definite, and its variational analysis
cannot be carried out by directly adapting the standard
fractional arguments. These features give rise to substantial
difficulties in each of the three regimes: mass-subcritical,
mass-critical, and mass-supercritical.

To overcome these difficulties, we rely on the estimates established
in \cite{Chen-Chen-Hauer} (see also \cite{BahrouniGouasmiaHajaiejOuannasser2026}), which provide essential control over the negative part of the logarithmic kernel and establish a quantitative link between the standard fractional Gagliardo seminorm and the positive logarithmic energy $[u]_{s+\log,+}^{2}$ (see Section~\ref{sec2}).
Precisely, the elementary bounds proved below are
\[
\mathcal J_-(u,u)\leq\frac{C(N,s)\sigma_{N-1}}{s^2}\|u\|_2^2,
\qquad
\mathcal J_s(u,u)\leq\frac{\mathcal J_+(u,u)}{C(N,s)\ln(1/r)}
+\frac{2\sigma_{N-1}}{sr^{2s}}\|u\|_2^2.
\]
Here $0<r<1$ and $\sigma_{N-1}=|\mathbb S^{N-1}|$. In particular, the squared fractional seminorm is infinitesimally controlled by the positive logarithmic form. This gives global constrained minimizers at every mass even at $p=2+4s/N$. To exclude vanishing, we use real oscillatory trial functions with Fourier support near the minimizing sphere of the symbol; they give energy strictly below the negative linear spectral bottom. Negativity alone would not provide this compactness argument.

A particularly delicate issue arises in the mass-supercritical regime. As in the fractional setting (see, for instance, \cite{LuoHajaiej2022,CherguiGouHajaiej2023,Soave2020,JeanjeanLuo2013} and the references therein), the energy functional is unbounded from below on the prescribed $L^2$-sphere, thereby precluding direct global minimization. Although constrained minimization on a Poho\u{z}aev manifold has been successfully employed in the aforementioned works, the logarithmic dependence of the energy on mass-preserving dilations gives rise to a substantially different fiber structure. More precisely, the presence of the term $t^{2s}\ln t$ significantly modifies the geometry of the energy along dilation paths, rendering the identification and analysis of the relevant Poho\u{z}aev branch considerably more delicate. A distinctive feature of our approach is the restriction of the variational analysis to the upper Poho\u{z}aev branch. Our strategy combines the exact fiber geometry with a uniform positive barrier at small prescribed norm. This motivates the introduction of the following constrained minimization problem:
\[
\gamma_p(c):=\inf_{u\in\mathcal{P}_c}\mathcal{E}_p(u),
\qquad c>0,
\]
where
\[
\mathcal{P}_c
:=
\left\{
u\in\mathbf{S}_c:\mathcal{Q}(u)=0
\right\},
\]
and
\[
\mathbf{S}_c
:=
\left\{
u\in\mathcal{S}_c:
2sC(N,s)\mathcal{J}_s(u,u)
+\frac{N(p-2)}{4p}
[4s-N(p-2)]
\|u\|_{L^p(\mathbb{R}^N)}^p<0
\right\}.
\]
Here, $\mathcal{Q}$ denotes the derivative of the energy functional with respect to the mass-preserving dilation, evaluated at $t=1$.
The strict inequality defining $\mathbf{S}_c$, together with $\mathcal Q(u)=0$, selects strict local fiber maxima. At a general point outside $\mathcal Q^{-1}(0)$, membership in $\mathbf S_c$ alone is not a critical-point assertion.

Without a positive fiber barrier, the remainder identity yields a one-sided characterization on $\{u\in\mathbf S_c:\mathcal Q(u)\geq0\}$ with $t\geq1$; local positivity to the left of $t=1$ does not justify a comparison along the entire orbit. At sufficiently small prescribed norm, however, a uniform positive barrier ensures two distinct Poho\u{z}aevscales on every fiber. The upper scale is then the positive global fiber maximum, and the full-fiber reduction in \eqref{globalreduced} becomes available.

The main ingredients of the small-mass existence proof are this uniform barrier, a smooth reduction to the upper branch, a strict comparison of the reduced levels as the prescribed norm varies, and compactness with exact recovery of mass. Together they give positive-energy upper-branch minimizers and compactness of every minimizing sequence modulo translations. The fiber-remainder identity also describes the geometry when the small-mass barrier is unavailable.
These results characterize small-mass upper-branch minimizers. For $N\geq2$, a second construction adapts the constrained symmetric linking method of \cite[Lemmas 2.3--2.5]{BartschValeriola} to the signed logarithmic form. Its operator-specific ingredients are a Jensen estimate ensuring boundedness at every mass, a high-energy multiplier bound restoring coercivity above the negative spectral bottom, and radial compactness up to the finite Sobolev endpoint. At that endpoint, the multiplier bound uses the fractional Sobolev inequality, since the subcritical estimate degenerates.

The small-mass limit has an additional feature tied to the nonmonotone symbol. Theorem~\ref{smallmass-spectral} shows that global minimizers in the fractional subcritical and critical ranges concentrate in Fourier space near $|\xi|=e^{-1/(2s)}$, while their mass-normalized profiles spread out in physical space. A weighted Fourier-shell interpolation estimate sharpens the energy, multiplier, and local-mass bounds. Anisotropic packets improve the lower gap bounds when $N>1$, and the Fourier decomposition also gives quantitative $L^q$ vanishing at the finite Sobolev endpoint. Corollary~\ref{one-dimensional-gaps} shows that the upper and lower energy-gap and multiplier-gap powers match in dimension one. Corollary~\ref{smallmass-physical} states the corresponding limits for the physical energy and stationary equation. These conclusions concern global minimizers in the subcritical and critical ranges, rather than the supercritical upper-branch minimizers. Corollary~\ref{green-kinetic-moment} also gives an exact scaling identity and the quantitative convergence of the normalized fractional kinetic moment to $1/e$.

Proposition~\ref{energy-mass-slopes} proves strict concavity of the minimum energy as a function of the actual mass $c^2$ and identifies its one-sided derivatives with the extreme minimizing multipliers. In particular, outside an at most countable set of prescribed norms, all global minimizers share the same multiplier, even though uniqueness of the profiles is not established. Corollary~\ref{multiplier-ordering} strengthens this comparison across different norms: every minimizing multiplier at a larger norm strictly exceeds every minimizing multiplier at a smaller norm. Thus a fixed multiplier determines at most one norm within this minimizing class.

\medskip
\noindent\textbf{Outline of the paper.} The paper is organized as follows: \textbf{Section 2} presents the preliminary results, introduces the functional framework, and states the main theorems. \textbf{Section 3} gives the common subcritical and critical existence argument and studies the small-mass energy gaps, spectral concentration, and spatial spreading. \textbf{Section 4} records the critical conclusion and explains the absence of a critical-mass obstruction. \textbf{Section 5} investigates the mass-supercritical regime on the upper Poho\u{z}aev branch through a fiber barrier argument, strict mass comparison, compactness modulo translations, and exact mass recovery. Moreover, for $N\geq2$, a symmetric radial minimax construction yields infinitely many normalized solutions at every prescribed norm throughout $2+4s/N<p\leq2_s^*$, with diverging energies. Finally, the \textbf{Appendix} collects auxiliary results and technical proofs used throughout the paper.

\section{Preliminaries and Main Results}\label{sec2}

In this section, we introduce the notation and preliminary results. Throughout, $N\geq1$, $0<s<1$, and $2<p<2_s^*$ unless stated otherwise; functions in the variational problems are real-valued. We write $\omega_N=|B_1|$ for the volume of the unit ball and $\sigma_{N-1}=|\mathbb S^{N-1}|=N\omega_N$ for the surface measure of the unit sphere. The Sobolev exponent is
\[
2_s^*=\begin{cases}2N/(N-2s),&N>2s,\\+\infty,&N\leq2s.\end{cases}
\]
When $2_s^*=\infty$, every exponent below it is understood to be finite.

The fractional Sobolev space $H^s(\mathbb R^N)$ consists of measurable real-valued functions and is defined by

\begin{equation*}
H^{s}(\mathbb{R}^{N})
:= \left\{
u \in L^{2}(\mathbb{R}^{N})
\,:\,
\iint_{\mathbb{R}^{2N}}
\frac{|u(x)-u(y)|^{2}}{|x-y|^{N+2s}}\,dx\,dy < \infty
\right\}.
\end{equation*}
This space is endowed with the norm
\begin{equation*}
\|u\|_{H^{s}(\mathbb{R}^{N})}^{2}
:= \|u\|_{L^{2}(\mathbb{R}^{N})}^{2} + [u]_{s}^{2},
\end{equation*}
where the Gagliardo seminorm is defined by
\begin{equation*}
[u]_{s}^{2}
:= \iint_{\mathbb{R}^{2N}}
\frac{|u(x)-u(y)|^{2}}{|x-y|^{N+2s}}\,dx\,dy.
\end{equation*}
Following \cite{Chen-Chen-Hauer, BahrouniGouasmiaHajaiejOuannasser2026}, we introduce the fractional logarithmic Sobolev space $H^{s+\log}(\mathbb{R}^{N})$, which is the energy space associated with the positive logarithmic kernel $\mathbf{k}^{s+\log,+}$. This kernel is defined by
\begin{equation*}
\mathbf{k}^{s+\log,+}(r)
:= \frac{C(N,s)\,(-\ln |r|)_{+}}{|r|^{N+2s}},
\qquad r \in \mathbb{R}^{N}\setminus\{0\},
\end{equation*}
where $a_+:=\max\{a,0\}$ and $a_-:=\max\{-a,0\}$. We also set $\mathbf k^{s+\log,-}(z)=C(N,s)(\ln|z|)_+/|z|^{N+2s}$. The space $H^{s+\log}(\mathbb{R}^{N})$ is defined as
\begin{equation*}
H^{s+\log}(\mathbb{R}^{N}) := \left\{ u \in L^{2}(\mathbb{R}^{N}) \; : \; [u]_{s+\log,+} < \infty \right\},
\end{equation*}
and it is endowed with the norm
\begin{equation*}
\|u\|_{H^{s+\log}(\mathbb{R}^{N})}^{2} := \|u\|_{L^{2}(\mathbb{R}^{N})}^{2} + [u]_{s+\log,+}^{2}.
\end{equation*}
Here, the associated Gagliardo-type seminorm is given by
\begin{equation*}
[u]_{s+\log,+}^{2} := \iint_{\mathbb{R}^{2N}} (u(x)-u(y))^{2}\, \mathbf{k}^{s+\log,+}(x-y)\, dx\, dy.
\end{equation*}
We recall that the space $H^{s+\log}(\mathbb{R}^{N})$ is a Hilbert space; see \cite[Proposition 3.1]{Chen-Chen-Hauer}. Moreover, it is separable and reflexive; see \cite[Proposition 2.3]{BahrouniGouasmiaHajaiejOuannasser2026}. For further details on additional properties and technical aspects of these spaces, we refer the reader to \cite{Chen-Chen-Hauer, BahrouniGouasmiaHajaiejOuannasser2026}. 
The quadratic-form energy convention used in the calculations below is defined by

\[
\mathcal{E}_{p}(u)
= \mathcal{E}_{s+\log}(u,u)
- \frac{1}{p}\int_{\mathbb{R}^{N}} |u|^{p}\,dx,
\]
and is considered under the constraint
\[
\mathcal{S}_{c}
= \left\{
u \in H^{s+\log}(\mathbb{R}^{N}) :
\int_{\mathbb{R}^{N}} |u|^{2}\,dx = c^{2}
\right\},
\qquad c>0.
\]
Here,
\[
\mathcal{E}_{s+\log}(u,u)
:= \mathcal{J}_{+}(u,u)-\mathcal{J}_{-}(u,u)
+ \frac{B(N,s)\,C(N,s)}{2}\,\mathcal{J}_{s}(u,u).
\]
The diagonal values of the bilinear forms $\mathcal{J}_{\pm}$ are defined by

\begin{equation*}
\mathcal{J}_{\pm}(u,u)
= \iint_{\mathbb{R}^{2N}}
(u(x)-u(y))^{2}\,
\mathbf{k}^{s+\log,\pm}(x-y)\,dx\,dy\newadd{.}
\end{equation*}
The corresponding bilinear values are obtained by replacing $(u(x)-u(y))^2$ by $(u(x)-u(y))(v(x)-v(y))$. The form $\mathcal{J}_{s}$ is given on the diagonal by

\begin{equation*}
\mathcal{J}_{s}(u,u)
= \iint_{\mathbb{R}^{2N}}
\frac{(u(x)-u(y))^{2}}{|x-y|^{N+2s}}\,dx\,dy.
\end{equation*}

\begin{remark}
	The full energy functional $\mathcal{E}_{s+\log}$ is not necessarily positive definite, which presents significant challenges in applying various tools from functional analysis. As established in \cite{Chen-Chen-Hauer, BahrouniGouasmiaHajaiejOuannasser2026}, for problems posed on bounded domains, the sign of the energy can be controlled by imposing suitable conditions on the domain and the constant $B(N,s)$. However, this approach is not applicable to the present problem, which is considered in the whole space $\mathbb{R}^{N}$.
\end{remark}

The logarithmic seminorm strengthens the fractional seminorm enough to restore local compactness at the finite Sobolev endpoint.
\begin{proposition}\label{proposition-embeddings}
For $0<s<1$ and $0<\epsilon<1-s$, the inclusions
\[
H^{s+\epsilon}(\R^N)\subset\HH\subset H^s(\R^N)
\]
are continuous. The embedding into $L^q_{\rm loc}$ is compact for $2\leq q\leq2_s^*$ if $N>2s$, and for every finite $q\geq2$ if $N\leq2s$.
\end{proposition}
\begin{proof}
The first inclusion follows from $(-\ln r)r^{2\epsilon}\leq C_\epsilon$ on $(0,1)$. For the second, split the fractional difference integral at $0<r<1$:
\begin{equation}\label{equ19}
\Js(w,w)\leq\frac{\Jp(w,w)}{C(N,s)\ln(1/r)}
       +\frac{2\sigma_{N-1}}{sr^{2s}}\|w\|_2^2.
\end{equation}
Thus, for every $\delta>0$, there is $D_\delta<\infty$ such that
\begin{equation}\label{local-infinitesimal}
\Js(w,w)\leq\delta\Jp(w,w)+D_\delta\|w\|_2^2.
\end{equation}
For endpoint compactness, take $u_n\rightharpoonup u$ in $\HH$, a compact set $K$, and a fixed $\eta\in C_c^\infty$ equal to one near $K$. Multiplication by $\eta$ is bounded on $\HH$: expand the difference of $\eta w$ and use
\[
\sup_y\int |\eta(x)-\eta(y)|^2\mathbf k^{s+\log,+}(x-y)\,dx<\infty.
\]
This integral is finite since $s<1$. Hence $w_n=\eta(u_n-u)$ is bounded in $\HH$, while local Rellich compactness gives $\|w_n\|_2\to0$. Applying \eqref{local-infinitesimal}, first letting $n\to\infty$ and then $\delta\downarrow0$, gives $\Js(w_n,w_n)\to0$. The fractional Sobolev inequality now implies $\|w_n\|_{2_s^*}\to0$ when $N>2s$. Interpolation proves the remaining exponents. If $N\leq2s$, the usual fractional Rellich theorem gives compactness into every finite $L^q_{\rm loc}$.
\end{proof}

\paragraph{Energy normalization.}
With the bilinear form defined above, $A(u):=\mathcal E_{s+\log}(u,u)=\langle(-\Delta)^{s+\log}u,u\rangle$. Thus the energy used in the variational calculations is $\E(u)=A(u)-\|u\|_p^p/p$, and its constrained critical points solve
\begin{equation}\tag{E}\label{variational-equation}
 2(-\Delta)^{s+\log}u+\lambda u=|u|^{p-2}u.
\end{equation}
The physical stationary equation \eqref{equ3} instead has energy $I(u)=A(u)/2-\|u\|_p^p/p$. The two conventions are related by $v=2^{-1/(p-2)}u$ and $\mu=\lambda/2$: a solution of \eqref{variational-equation} of norm $c$ gives a solution of \eqref{equ3} of norm $2^{-1/(p-2)}c$. 
We retain $\E$ in the calculations below and use this mass conversion in all three regimes. In particular, to solve the physical problem at norm $c$, apply a theorem for $\E$ at norm $2^{1/(p-2)}c$. Writing $k_p=2^{-1/(p-2)}$, the energy conversion is
\begin{equation}\label{green-normalization}
 I(k_pu)=\frac{k_p^2}{2}\mathcal E_p(u),\qquad
 \|k_pu\|_2=k_p\|u\|_2,\qquad \mu=\frac\lambda2.
\end{equation}
Thus the conversion preserves constrained minimization, compactness, and the sign of the energy, while rescaling the prescribed norm.}
The functional $\E$ is of class $C^1$ on $\HH$, so its restriction to $\Sc$ is well-defined. We consider
\begin{equation}\label{equ}
\mathcal{M}_{p}(c)
:= \inf_{u \in \mathcal{S}_{c}} \mathcal{E}_{p}(u),
\qquad c>0.
\end{equation}

A weak solution of \eqref{equ3} means a real $u\in\HH$ satisfying
\[
\mathcal E_{s+\log}(u,\varphi)+\lambda\int u\varphi
=\int |u|^{p-2}u\varphi\qquad(\varphi\in\HH).
\]
For \eqref{variational-equation}, the first bilinear term is multiplied by two. All terms are finite by the embeddings above. We recall the ground-state definition used in \cite{BellazziniJeanjean2016}.

\begin{definition} \label{definition1}
	A function $u_{c} \in \mathcal{S}_{c}$ is called a ground-state solution of Problem \eqref{variational-equation} if it attains the minimum of the energy functional among all constrained critical points in $\mathcal{S}_{c}$. More precisely,
	\[
	\mathcal{E}_{p}(u_{c})
	=
	\inf\left\{
	\mathcal{E}_{p}(v):
	v \in \mathcal{S}_{c},\;
	\left(\mathcal{E}_{p}\big|_{\mathcal{S}_{c}}\right)'(v)=0
	\right\}.
	\]
	\end{definition}

\begin{remark}
	This definition also applies when $\mathcal E_p$ is unbounded below on $\Sc$. We also note that every minimizer of $\mathcal{M}_{p}(c)$ is a ground-state solution of Problem~\eqref{variational-equation}.
\end{remark}

For $0<s<1$ and $2\leq p<2_s^*$, the fractional logarithmic Gagliardo--Nirenberg inequality follows from the usual fractional inequality \cite{Hajaiej2011} and the continuous inclusion $\HH\subset H^s$. For $N>2s$, it also follows by setting the base exponent equal to two in \cite[Theorem 3.6]{BahrouniGouasmiaHajaiejOuannasser2026}. Thus there exists $\mathbf K=\mathbf K(N,s,p)>0$ such that, for every $u\in\HH$,

\begin{equation}\label{equ13}
\|u\|_{L^{p}(\mathbb{R}^{N})}^{p}
\leq
\mathbf{K}\,
\|u\|_{H^{s+\log}(\mathbb{R}^{N})}^{\frac{N(p-2)}{2s}}
\|u\|_{L^{2}(\mathbb{R}^{N})}^{p-\frac{N(p-2)}{2s}}.
\end{equation}
At the finite endpoint $p=2_s^*$, the same inequality holds with zero $L^2$ exponent, and its sharp constant is attained by Proposition~\ref{proposition1}. We denote by $\mathbf{K}_{\mathrm{opt}}$ the sharp constant in this inequality, which depends on $s$, $N$, and $p$ and is obtained by minimizing the associated ``Weinstein functional'' (see Proposition~\ref{proposition1}). Finally, we shall also invoke the following result:
\begin{lemma}\label{lemma}
Let $2<p<2_s^*$ and $u_n\rightharpoonup u$ in $\HH$. Then
\begin{equation}\label{equA1}
\E(u_n)=\E(u_n-u)+\E(u)+o(1).
\end{equation}
The same splitting holds for each of $\Jp$, $\Jm$, and $\Js$.
\end{lemma}
\begin{proof}
The bilinear form associated with each quadratic form is continuous on $\HH$. The form $\Jm$ is continuous even on $L^2$. Polarization gives, with $v_n=u_n-u$,
\[
\mathcal J_\pm(u_n,u_n)-\mathcal J_\pm(v_n,v_n)-\mathcal J_\pm(u,u)
=2\mathcal J_\pm(v_n,u)\to0.
\]
The argument for $\Js$ is identical. Local compactness gives almost-everywhere convergence along a subsequence, so the Br\'ezis--Lieb lemma gives the $L^p$ splitting. The conclusion holds for the full sequence: any subsequence has a further subsequence with this property. Combining the splittings proves \eqref{equA1}.
\end{proof}

Our first main result in the mass-subcritical case is stated as follows.

\begin{theorem}\label{theorem1}
	Assume that $N \geq 1$, $0 < s < 1$, and \(2 < p < 2+\frac{4s}{N}\). Then, for every $c>0$, the infimum $\mathcal M_p(c)$ is attained by a function $u_c$. Moreover, every minimizing sequence is relatively compact in $\HH$ modulo translations, and there exists a Lagrange multiplier \(\lambda_c>2/(es)\) such that \(u_c\) is a ground-state solution of \eqref{variational-equation}--\eqref{equ5} with \(\lambda=\lambda_c\).
\end{theorem}

The second main result in the mass-critical case is stated as follows.

\begin{theorem}\label{theorem3}
Assume $N\geq1$, $0<s<1$, and $p=2+4s/N$. For every $c>0$, the infimum $\mathcal M_p(c)$ is attained by a real function $u_c\in\Sc$. Every minimizing sequence is relatively compact in $\HH$ modulo translations. The minimizer is a ground state for \eqref{variational-equation}, with multiplier $\lambda_c>2/(es)$. Consequently, the physical equation \eqref{equ3} also has a normalized global energy minimizer at every prescribed norm, after the amplitude conversion above.
\end{theorem}

Theorems~\ref{theorem1} and \ref{theorem3} use a common compactness argument. The crucial comparison is $\mathcal M_p(c)<-c^2/(es)$, rather than merely $\mathcal M_p(c)<0$. At the critical power, \eqref{local-infinitesimal} absorbs the nonlinear term at every fixed mass, so no critical-mass threshold is necessary.

In the mass-supercritical case, direct minimization in \eqref{equ} is unavailable because the constrained energy is unbounded below. Indeed, for any $u \in \mathcal{S}_{c}$ and $t>0$, we define
\[
u_{t}(x):=t^{\frac{N}{2}}u(tx), \qquad t>0.
\]
A straightforward calculation shows that $u_{t}\in\mathcal{S}_{c}$ for every $t>0$. Moreover,
\begin{equation*}
\begin{aligned}
\mathcal{E}_{p}(u_{t})
&= \mathcal{E}_{s+\log}(u_{t},u_{t})
-\frac{1}{p}\int_{\mathbb{R}^{N}}|u_{t}|^{p}\,dx\\[4pt]
&= \mathcal{J}_{+}(u_{t},u_{t})
-\mathcal{J}_{-}(u_{t},u_{t})
+\frac{B(N,s)\,C(N,s)}{2}\,\mathcal{J}_{s}(u_{t},u_{t})
-\frac{1}{p}\int_{\mathbb{R}^{N}}|u_{t}|^{p}\,dx\\[4pt]
& =
t^{2s}
\left[
\mathcal{J}_{+}(u,u)
-\mathcal{J}_{-}(u,u)
+
C(N,s)
\left(
\ln t+\frac{B(N,s)}{2}
\right)
\mathcal{J}_{s}(u,u)
\right]
\\
&\quad
-\frac{t^{N\left(\frac{p}{2}-1\right)}}{p}
\int_{\mathbb{R}^{N}}|u|^{p}\,dx.
\end{aligned}
\end{equation*}
Since $ p > 2+\frac{4s}{N}, $ we consequently obtain
\[
\mathcal{E}_{p}(u_{t})\to-\infty
\qquad\text{as }t\to\infty.
\]
It follows that
\[
\mathcal{M}_{p}(c)=-\infty
\qquad\text{for every }c>0.
\]
Therefore, to establish the existence of solutions to
\eqref{variational-equation}--\eqref{equ5}, we introduce the following minimization
problem:
\begin{equation}\label{equu}
\gamma_{p}(c)
:=\inf_{u\in\mathcal{P}_{c}}\mathcal{E}_{p}(u),
\qquad c>0,
\end{equation}
where $\mathcal P_c$ denotes the upper Poho\u{z}aev branch, defined by
\begin{equation}\label{equ27}
\mathcal{P}_{c}
:=\left\{
u\in\textbf{S}_c:\mathcal{Q}(u)=0
\right\},
\end{equation}
with
\[
\textbf{S}_c := \left\lbrace u \in \mathcal{S}_c:\quad
2s\, C(N,s)\mathcal{J}_{s}(u,u)
+\frac{N(p-2)}{4p}
\left[4s-N(p-2)\right]
\|u\|_{L^{p}(\mathbb{R}^{N})}^{p} < 0
\right\rbrace,
\]
and where the functional $\mathcal{Q}$ is defined as the derivative
of $\mathcal{E}_{p}(u_t)$ with respect to $t$, evaluated at $t=1$:
\[
\begin{aligned}
\mathcal{Q}(u)
&:=\left.\frac{d}{dt}\mathcal{E}_{p}(u_{t})\right|_{t=1}\\[4pt]
&=
2s\left[
\mathcal{J}_{+}(u,u)-\mathcal{J}_{-}(u,u)
\right]
+C(N,s)\bigl(1+sB(N,s)\bigr)\mathcal{J}_{s}(u,u)\\[4pt]
&\quad
-\frac{N(p-2)}{2p}
\int_{\mathbb{R}^{N}}|u|^{p}\,dx.
\end{aligned}
\]

The constraint $\mathcal Q(u)=0$ is the scaling identity for the variational energy. Its validity for the minimizing profiles is proved directly in Lemma~\ref{profilevariation}; the role of admissible dilation tests is recorded in Lemma~\ref{lemma2}.

   Our next main result is stated as follows.

\begin{theorem}[Small-mass upper-branch minimization]\label{theorem4}
Assume that $N\geq1$, $0<s<1$, and $2+4s/N<p<2_s^*$, and let $c_*$ be the positive threshold in \eqref{newconstants}. If $0<c<c_*$, then:
\begin{enumerate}
\item Every $u\in\Sc$ has exactly two Poho\u{z}aevscales. At its unique upper-branch scale $t_+(u)$, the projected function $u_{t_+(u)}$ has positive energy. In particular, $\Pc\ne\varnothing$.
\item The infimum $\gamma_p(c)$ is positive and is attained by a function $u_c\in\Pc$. Every minimizing sequence is relatively compact in $\HH$ after translations and passage to a subsequence. In particular, the limit has exactly the prescribed $L^2$-norm $c$.
\item The minimizer is a constrained critical point of $\E$ on $\Sc$ and satisfies the Euler--Lagrange equation below, where $\Lambda$ is defined in \eqref{newconstants}:
	\begin{equation}\label{equationcorrect}
	2 (-\Delta)^{s+\log} u_c+\lambda_c u_c=|u_c|^{p-2}u_c\quad\text{in }\R^N,
	\qquad \lambda_c>\Lambda>\frac2{es}.
	\end{equation}
\item The map $m\mapsto\gamma_p(m)$ is strictly decreasing on $(0,c_*)$. For $0<m_1<m_2<c_*$, we have
	\begin{equation}\label{strictquant}
	\gamma_p(m_2)\leq\gamma_p(m_1)
	-\frac{C}{s}\int_{m_1}^{m_2}\frac{\rho_m}{m}\,dm
	<\gamma_p(m_1),
	\end{equation}
where $\rho_m$ is defined in \eqref{newconstants0} and $C=C(N,s)$.

\end{enumerate}
The theorem makes no assertion about the sign, radial symmetry, or uniqueness of the minimizer, or about its energy relative to all normalized critical points.
\end{theorem}

\begin{theorem}[Radial normalized solutions at every mass]\label{theorem-all-masses}
Assume that $N\geq2$, $0<s<1$, and $2+\frac{4s}{N}<p\leq2_s^*.$ The finite Sobolev-critical endpoint is included. For every prescribed $L^2$-norm $c>0$, there exist infinitely many distinct pairs of real radial weak solutions $\{u_{c,j},-u_{c,j}\}\subset H^{s+\log}(\R^N)$ and real numbers $\lambda_{c,j}$ such that
\[
 (-\Delta)^{s+\log}u_{c,j}+\lambda_{c,j}u_{c,j}
       =|u_{c,j}|^{p-2}u_{c,j}\quad\text{in }\R^N,
 \qquad \|u_{c,j}\|_2=c.
\]
With the physical energy $I(u)=\mathcal E_{s+\log}(u,u)/2-\|u\|_p^p/p$, these solutions can be chosen so that
\[
 I(u_{c,j})\longrightarrow+\infty,\qquad
 \lambda_{c,j}\longrightarrow+\infty.
\]
In particular, infinitely many have positive energy and satisfy $\lambda_{c,j}>1/(es)$. No assertion that these functions are positive or that they minimize the upper-branch energy is made.
\end{theorem}
The small-mass minimization theorem and the all-mass existence theorem have different conclusions. The latter is obtained from high radial minimax levels and does not require positivity of $\gamma_p(c)$.

\begin{corollary}[Physical normalized minimizers]\label{review-physical}
For $N\geq1$, $0<s<1$, $2<p\leq2+4s/N$, and every $c>0$, the physical energy
\[
I(v)=\tfrac12\langle(-\Delta)^{s+\log}v,v\rangle-\tfrac1p\|v\|_p^p
\]
attains its minimum on $\|v\|_2=c$. Every minimizing sequence is relatively compact in $\HH$ modulo translations, and each minimizer solves \eqref{equ3} with $\lambda>1/(es)$.
\end{corollary}
\begin{proof}
Set $k_p=2^{-1/(p-2)}$. The map $u\mapsto k_pu$ is a bijection between the spheres of norms $c/k_p$ and $c$, and
\[
I(k_pu)=\tfrac{k_p^2}{2}\E(u).
\]
Apply Theorems~\ref{theorem1} and \ref{theorem3}, together with the sharper multiplier estimate \eqref{equ14}.
\end{proof}

\section{Global Minimizers and Small-Mass Asymptotics}

We prove compactness and existence for the quadratic-form energy $\E$. The physical energy is obtained by the amplitude conversion stated above. Put
\[
A(u)=\mathcal E_{s+\log}(u,u),\quad K(u)=\Jp(u,u),\quad
S(u)=\Js(u,u),\quad P(u)=\|u\|_p^p.
\]
With the unitary Fourier transform,
\begin{equation}\label{review-spectral}
A(u)=\int\ell(\xi)|\widehat u(\xi)|^2\,d\xi,\qquad
\ell(\xi)=|\xi|^{2s}\ln|\xi|^2\geq-\frac1{es}.
\end{equation}
The negative logarithmic tail satisfies
\begin{equation}\label{equ11}
\Jm(u,u)\leq \frac{C(N,s)\sigma_{N-1}}{s^2}\|u\|_2^2.
\end{equation}
Equations \eqref{equ19} and \eqref{equ11} imply
\begin{equation}\label{review-Acoercive}
A(u)\geq\tfrac12K(u)-D\|u\|_2^2
\end{equation}
for a constant $D$. More generally, for every $L>0$, $A(u)\geq LS(u)-D_L\|u\|_2^2$.

\begin{lemma}\label{lemma-coercivity}
If $2<p<2+4s/N$, then $\E$ is coercive on $\Sc$ for every $c>0$. The same conclusion holds when $p=2+4s/N$.
\end{lemma}
\begin{proof}
The fractional Gagliardo--Nirenberg inequality gives
\begin{equation}\label{review-GN}
P(u)\leq G S(u)^r c^\beta,\qquad
r=\frac{N(p-2)}{4s},\quad \beta=p-2r>0.
\end{equation}
If $r<1$, the continuous embedding estimate $S\leq C(K+c^2)$ and \eqref{review-Acoercive} yield
\begin{equation}\label{equ7}
\E(u)\geq\tfrac12K(u)-Dc^2-C_c(K(u)+c^2)^r\to+\infty
\end{equation}
as $K(u)\to\infty$. If $r=1$, use \eqref{local-infinitesimal} instead and choose $\delta>0$ with $G c^{p-2}\delta/p\leq1/4$. Then
\begin{equation}\label{review-critical-coercive}
\E(u)\geq\tfrac14K(u)-D_c,
\end{equation}
where $D_c<\infty$ for each fixed $c$. This proves the assertion at every mass, including the critical power.
\end{proof}

\begin{lemma}[Strict comparison with the spectral bottom]\label{review-trial}
If $2<p\leq2+4s/N$, then, for every $c>0$,
\begin{equation}\label{review-strict-bottom}
-\infty<\mathcal M_p(c)<-\frac{c^2}{es}.
\end{equation}
In particular, every minimizing sequence satisfies $\liminf P(u_n)>0$.
\end{lemma}
\begin{proof}
Choose a real Schwartz function $\phi$ with $\|\phi\|_2=1$ whose Fourier transform is even, smooth, and compactly supported. Fix $\xi_0\ne0$ with $|\xi_0|=e^{-1/(2s)}$, so $\ell(\xi_0)=-1/(es)$ and $\nabla\ell(\xi_0)=0$. For sufficiently large $R$, set
\[
u_R(x)=\sqrt2\,cR^{-N/2}\phi(x/R)\cos(\xi_0\cdot x).
\]
The two Fourier packets centered at $\pm\xi_0$ have disjoint supports, so $\|u_R\|_2=c$ exactly. Taylor expansion of $\ell$ near these centers gives
\[
A(u_R)=-\frac{c^2}{es}+O(R^{-2}).
\]
Writing $a=N(p-2)/2$, periodic averaging gives
\[
P(u_R)=(\sqrt2c)^pR^{-a}
\left(\mu_p\|\phi\|_p^p+o(1)\right),\qquad
\mu_p=\frac1\pi\int_0^\pi|\cos t|^p\,dt>0.
\]
For completeness, this averaging follows by uniformly approximating $|\cos t|^p$ by trigonometric polynomials and applying the Riemann--Lebesgue lemma to the nonconstant modes against $|\phi|^p\in L^1$. Since $a\leq2s<2$, the negative nonlinear term dominates the $O(R^{-2})$ spectral error. Thus $\E(u_R)<-c^2/(es)$ for large $R$. Coercivity gives a finite lower bound for the infimum. Finally, $A(u_n)\geq-c^2/(es)$ implies
\[
\liminf P(u_n)\geq p\left(-\frac{c^2}{es}-\mathcal M_p(c)\right)>0.
\]
\end{proof}

\begin{lemma}\label{lemma-mass-continuity}
For $2<p\leq2+4s/N$, the map $c\mapsto\mathcal M_p(c)$ is continuous and strictly decreasing on $(0,\infty)$.
\end{lemma}
\begin{proof}
For $t>1$, amplitude multiplication gives
\[
\mathcal M_p(tc)\leq t^2\mathcal M_p(c)
-\frac{t^p-t^2}{p}\liminf P(u_n)<t^2\mathcal M_p(c).
\]
Here $(u_n)$ is a minimizing sequence at norm $c$, and Lemma~\ref{review-trial} makes the inequality strict. Since $\mathcal M_p(c)<0$, strict decrease follows. To prove continuity, let $c_n\to c>0$. A fixed near-minimizer at $c$, multiplied by $c_n/c$, gives $\limsup\mathcal M_p(c_n)\leq\mathcal M_p(c)$. Conversely, choose near-minimizers $v_n\in\mathcal S_{c_n}$. The coercivity estimates are uniform when $c_n$ lies in a compact subinterval of $(0,\infty)$, so $(v_n)$ is bounded in $\HH$. Therefore
\[
\mathcal M_p(c)\leq\E((c/c_n)v_n)=\mathcal M_p(c_n)+o(1).
\]
This gives the reverse inequality and continuity.
\end{proof}

\begin{lemma}[Strict subadditivity in squared mass]\label{lemma1}
Let $2<p\leq2+4s/N$. For $0<a<c$ and $b=(c^2-a^2)^{1/2}$,
\begin{equation}\label{review-binding}
\mathcal M_p(c)<\mathcal M_p(a)+\mathcal M_p(b).
\end{equation}
\end{lemma}
\begin{proof}
Lemma~\ref{lemma-mass-continuity} actually proves $\mathcal M_p(tc)<t^2\mathcal M_p(c)$ for $t>1$. Apply this with $t=c/a$ and $t=c/b$ to get
\[
\frac{a^2}{c^2}\mathcal M_p(c)<\mathcal M_p(a),\qquad
\frac{b^2}{c^2}\mathcal M_p(c)<\mathcal M_p(b).
\]
Addition proves the result because $a^2+b^2=c^2$.
\end{proof}

\begin{proof}[Proof of Theorem~\ref{theorem1}]
The argument applies throughout the range $2<p\leq2+4s/N$. Let $(u_n)\subset\Sc$ minimize $\mathcal M_p(c)$. Coercivity gives boundedness in $\HH$. If
\[
\sup_y\int_{B_1(y)}|u_n|^2\to0,
\]
the fractional Lions vanishing criterion (cf.~\cite{Lions}) gives $P(u_n)\to0$, contradicting Lemma~\ref{review-trial}. Here the criterion applies to bounded $H^s$ sequences and exponents $2<p<2_s^*$; it follows by local fractional interpolation on unit cubes, summing with bounded overlap, and interpolation with the global Sobolev bounds. When $N\leq2s$, use a finite Sobolev exponent larger than $p$. Thus, after translations and extraction,
\[
u_n\rightharpoonup u\ne0\text{ in }\HH,\qquad
u_n\to u\text{ locally in }L^p\text{ and a.e.}
\]
Put $a=\|u\|_2$. If $a<c$, Hilbert splitting gives $\|u_n-u\|_2\to b=(c^2-a^2)^{1/2}>0$. By energy splitting and continuity of the levels,
\[
\mathcal M_p(c)=\E(u)+\E(u_n-u)+o(1)
\geq\mathcal M_p(a)+\mathcal M_p(b)+o(1),
\]
contradicting \eqref{review-binding}. Hence $a=c$ and $v_n=u_n-u\to0$ in $L^2$. Boundedness and \eqref{local-infinitesimal} imply $S(v_n)\to0$, while \eqref{equ11} gives $\Jm(v_n,v_n)\to0$. Interpolation yields $P(v_n)\to0$. Consequently energy splitting becomes
\[
\mathcal M_p(c)=\E(u)+K(v_n)+o(1).
\]
Since $\E(u)\geq\mathcal M_p(c)$ and $K(v_n)\geq0$, it follows that $K(v_n)\to0$ and $\E(u)=\mathcal M_p(c)$. Thus convergence is strong in $\HH$, proving compactness modulo translations for every minimizing sequence.

Constrained criticality gives \eqref{variational-equation}, and testing with $u$ gives
\begin{equation}\label{equ160}
2A(u)+\lambda_c c^2=P(u).
\end{equation}
In fact,
\begin{equation}\label{equ14}
\lambda_c c^2=\left(1-\frac2p\right)P(u)-2\mathcal M_p(c)
>\frac{2c^2}{es},
\end{equation}
so \newadd{$\lambda_c>2/(es)$}. A global constrained minimizer is a ground state by Definition~\ref{definition1}. Finally, the amplitude map $v=2^{-1/(p-2)}u$ gives the corresponding physical solution. Applying the argument at norm $2^{1/(p-2)}c$ gives a physical solution at any prescribed norm $c$.
\end{proof}

\subsection*{Dependence of the minimum energy on the prescribed norm}
\begin{proposition}[Energy slopes and minimizing multipliers]\label{energy-mass-slopes}
Let $N\geq1$, $0<s<1$, and $2<p\leq2+4s/N$. For $c>0$ define
\[
 \mathcal G_c=\{u\in\mathcal S_c:\mathcal E_p(u)=\mathcal M_p(c)\},
 \qquad
 \lambda(u)=\frac{P(u)-2A(u)}{c^2}\quad(u\in\mathcal G_c).
\]
In terms of the actual mass $m=c^2$, the minimum level
\[
 \mathfrak e(m):=\mathcal M_p(\sqrt m),\qquad m>0,
\]
is strictly concave. The extrema $\lambda_{\min}(c)$ and $\lambda_{\max}(c)$ of $\lambda$ on $\mathcal G_c$ are attained. The function $\mathcal M_p$ is locally Lipschitz on $(0,\infty)$, and its right and left derivatives satisfy
\begin{equation}\label{mass-one-sided-slopes}
 \mathcal M'_{p,+}(c)=-c\lambda_{\max}(c),\qquad
 \mathcal M'_{p,-}(c)=-c\lambda_{\min}(c).
\end{equation}
It is differentiable outside an at most countable set. At every differentiability point, all global minimizers have the same multiplier and
\begin{equation}\label{mass-slope-identity}
 \mathcal M_p'(c)=-c\lambda_c<-\frac{2c}{es}.
\end{equation}
In the physical convention, if $m_I(c)=\inf_{\|w\|_2=c}I(w)$, then at every differentiability point
\[
 m_I'(c)=-c\mu_c<-\frac{c}{es},
\]
where $\mu_c$ is the multiplier of any physical minimizer. These assertions concern the minimum level and multipliers, and do not imply uniqueness of minimizing profiles.
\end{proposition}
\begin{proof}
Amplitude normalization gives
\begin{equation}\label{concave-squared-mass}
 \mathfrak e(m)=\inf_{\|v\|_2=1}
 \left(mA(v)-\frac{m^{p/2}}pP(v)\right).
\end{equation}
For each fixed unit-norm $v$, the expression in parentheses is strictly concave in $m>0$, since $p/2>1$ and $P(v)>0$. For distinct $m_1,m_2>0$ and $0<\theta<1$, choose a minimizer at $m_\theta=\theta m_1+(1-\theta)m_2$ and normalize it to a unit-norm profile $v_\theta$. Strict concavity of its scalar expression gives
\[
 \begin{aligned}
 \mathfrak e(m_\theta)
 &>\theta\left(m_1A(v_\theta)-\frac{m_1^{p/2}}pP(v_\theta)\right)\\
 &\quad+(1-\theta)\left(m_2A(v_\theta)-\frac{m_2^{p/2}}pP(v_\theta)\right)\\
 &\geq\theta\mathfrak e(m_1)+(1-\theta)\mathfrak e(m_2).
 \end{aligned}
\]
Thus $\mathfrak e$ is finite and strictly concave. A finite concave function on an open interval is locally Lipschitz and has unequal one-sided derivatives at at most countably many points: each such point accounts for a distinct jump interval of its monotone derivatives. The identity $\mathcal M_p(c)=\mathfrak e(c^2)$ transfers these regularity conclusions to $\mathcal M_p$. For the comparisons below, we also record the amplitude-normalized level
\begin{equation}\label{concave-rescaled-level}
 F(t):=t^{-2/(p-2)}\mathcal M_p\bigl(t^{1/(p-2)}\bigr)
 =\inf_{\|v\|_2=1}\left(A(v)-\frac{t}{p}P(v)\right),\qquad t>0.
\end{equation}

Compactness modulo translations of $\mathcal G_c$, and continuity and translation invariance of $A$ and $P$, give attainment of the multiplier extrema. We also need compactness when the norm varies. If $c_n\to c>0$ and $w_n\in\mathcal G_{c_n}$, the uniform coercivity estimates on compact norm intervals bound $(w_n)$ in $\HH$. Continuity of $\mathcal M_p$ and amplitude scaling then show that $(c/c_n)w_n$ is a minimizing sequence at norm $c$. Theorems~\ref{theorem1} and \ref{theorem3} therefore give convergence, after translations and extraction, to a member of $\mathcal G_c$.

For $u\in\mathcal G_c$, the amplitude competitor gives
\[
 \mathcal M_p(c+h)-\mathcal M_p(c)
 \leq\mathcal E_p\bigl((1+h/c)u\bigr)-\mathcal E_p(u)
 =-c\lambda(u)h+O(h^2).
\]
For $h>0$, this bounds the right derivative above by $-c\lambda_{\max}(c)$. Conversely, take $w_h\in\mathcal G_{c+h}$ and compare $(c/(c+h))w_h$ with the minimum at norm $c$. Expansion gives
\[
 \frac{\mathcal M_p(c+h)-\mathcal M_p(c)}h
 \geq\frac{2A(w_h)-P(w_h)}{c+h}+O(h).
\]
The remainder is uniform for small $h$, by the preceding boundedness. Compactness at varying norm shows that every limiting right-hand side is $-c\lambda(u)$ for some $u\in\mathcal G_c$, and hence is at least $-c\lambda_{\max}(c)$. This proves the first formula in \eqref{mass-one-sided-slopes}. The same comparisons with $h<0$, reversing the inequalities upon division by $h$, prove the second. Equality of the derivatives is equivalent to equality of the multiplier extrema. The strict bound follows from \eqref{equ14}. In the mass variable, the same formulas read
\[
 \mathfrak e'_+(m)=-\frac12\lambda_{\max}(\sqrt m),\qquad
 \mathfrak e'_-(m)=-\frac12\lambda_{\min}(\sqrt m).
\]

Finally, with $k_p=2^{-1/(p-2)}$,
\[
 m_I(c)=\frac{k_p^2}{2}\mathcal M_p(c/k_p),\qquad
 \mu_c=\frac12\lambda_{c/k_p}.
\]
Differentiating at a differentiability point gives the stated physical identity.
\end{proof}

\begin{corollary}[Strict ordering of minimizing multipliers]\label{multiplier-ordering}
Under the assumptions of Proposition~\ref{energy-mass-slopes}, let $0<c_1<c_2$ and choose arbitrary $u_i\in\mathcal G_{c_i}$. Write $v_i=u_i/c_i$ and $\lambda_i=\lambda(u_i)$. Then
\begin{equation}\label{multiplier-ordering-bound}
 \lambda_2-\lambda_1
 \geq\bigl(c_2^{p-2}-c_1^{p-2}\bigr)\|v_1\|_p^p>0.
\end{equation}
In particular,
\[
 \lambda_{\max}(c_1)<\lambda_{\min}(c_2).
\]
Uniformly over all choices of global minimizers,
\[
 \lambda(u_c)\longrightarrow\frac2{es}\quad(c\downarrow0),
 \qquad
 \lambda(u_c)\longrightarrow+\infty\quad(c\to\infty).
\]
The same strict ordering holds for the physical minimizing multipliers $\mu_c$, whose small-norm limit is $1/(es)$ and whose large-norm limit is $+\infty$. Consequently, a fixed multiplier can correspond to at most one prescribed norm among global constrained minimizers. This conclusion does not assert uniqueness of profiles or extend to arbitrary normalized critical points.
\end{corollary}
\begin{proof}
Put $t_i=c_i^{p-2}$ and $P_i=\|v_i\|_p^p$. In \eqref{concave-rescaled-level}, $v_i$ minimizes $A(v)-t_iP(v)/p$ on $\|v\|_2=1$. Comparing the two minimizers in both problems gives
\[
 A(v_1)-\frac{t_1}{p}P_1\leq A(v_2)-\frac{t_1}{p}P_2,
 \qquad
 A(v_2)-\frac{t_2}{p}P_2\leq A(v_1)-\frac{t_2}{p}P_1.
\]
Adding yields $(t_2-t_1)(P_2-P_1)\geq0$, so $P_2\geq P_1>0$. The second comparison also yields
\[
 F(t_2)-F(t_1)\leq-\frac{t_2-t_1}{p}P_1.
\]
The multiplier identity is
\[
 \lambda_i=\left(1-\frac2p\right)t_iP_i-2F(t_i).
\]
Hence
\[
 \begin{aligned}
 \lambda_2-\lambda_1
 &\geq\left(1-\frac2p\right)(t_2P_2-t_1P_1)
       +\frac2p(t_2-t_1)P_1\\
 &\geq(t_2-t_1)P_1>0,
 \end{aligned}
\]
proving \eqref{multiplier-ordering-bound}. Taking minimizers attaining the multiplier extrema gives the asserted ordering of those extrema.

The small-norm limit is Theorem~\ref{smallmass-spectral}. To prove the large-norm limit, fix any $\phi\in\HH$ with $\|\phi\|_2=1$. By \eqref{concave-rescaled-level},
\[
 F(t)\leq A(\phi)-\frac{t}{p}P(\phi).
\]
Since $\lambda(u_c)=(1-2/p)c^{p-2}P(u_c/c)-2F(c^{p-2})$, it follows that
\[
 \lambda(u_c)\geq\frac2p c^{p-2}P(\phi)-2A(\phi)
 \longrightarrow+\infty.
\]
This bound is independent of the minimizing choice. The physical conclusions follow from the fixed norm conversion and $\mu_c=\lambda_{c/k_p}/2$. Strict ordering rules out two different minimizing norms at a common multiplier.
\end{proof}

\subsection*{The small-mass limit: a spectral shell and spatial spreading}

\begin{lemma}[Interpolation near the minimizing Fourier sphere]\label{shell-interpolation}
Let $2<p\leq2+4s/N$. There are constants $C$ and, for each $R>0$, $C_R$ such that every $v\in\HH$ with $\|v\|_2=1$ and
\[
 0<D(v):=\int_{\R^N}(\ell(\xi)+b_s)|\widehat v(\xi)|^2\,d\xi\leq1,
 \qquad b_s=\frac1{es},
\]
satisfies
\begin{equation}\label{shell-interpolation-bounds}
 \|v\|_p^p\leq C D(v)^{(p-2)/4},\qquad
 \sup_y\int_{B_R(y)}|v|^2\leq C_R D(v)^{1/2}.
\end{equation}
\end{lemma}
\begin{proof}
Put $h(\xi)=\ell(\xi)+b_s$ and $D=D(v)$. Fix a sufficiently thin annulus
$\mathcal A=\{\xi:||\xi|-r_s|<\epsilon_0\}$ around $r_s=e^{-1/(2s)}$. The radial symbol has
$\ell''(r_s)=4s r_s^{2s-2}>0$. Thus $h(\xi)\geq k_s(|\xi|-r_s)^2$ on $\mathcal A$, and $h\geq k_0>0$ outside it. Polar integration, also interpreted as the sum over the two intervals when $N=1$, gives
\begin{equation}\label{shell-resolvent-integral}
 \int_{\mathcal A}\frac{d\xi}{h(\xi)+D}
 \leq C_s\int_{-\epsilon_0}^{\epsilon_0}\frac{dt}{k_st^2+D}
 \leq C_sD^{-1/2}.
\end{equation}
Write $v=v_{\mathcal A}+v_{\mathcal A^c}$ using these Fourier restrictions. Weighted Cauchy--Schwarz and Fourier inversion yield
\[
 \|v_{\mathcal A}\|_\infty^2
 \leq C\left(\int_{\mathcal A}(h+D)|\widehat v|^2\right)
            \left(\int_{\mathcal A}\frac{d\xi}{h+D}\right)
 \leq C D^{1/2},
\]
because the first integral is at most $2D$. Also
$\|v_{\mathcal A^c}\|_2^2\leq D/k_0$. The coercive estimate \eqref{review-Acoercive}, applied to $v_{\mathcal A^c}$, gives
\[
 \|v_{\mathcal A^c}\|_{\HH}^2\leq CD.
\]
Indeed its signed quadratic energy is at most $D$, and its $L^2$ norm is already controlled by $D$. The continuous $L^p$ embedding and $\|v_{\mathcal A}\|_2\leq1$ imply
\[
 \|v_{\mathcal A}\|_p^p\leq C D^{(p-2)/4},\qquad
 \|v_{\mathcal A^c}\|_p^p\leq C D^{p/2}.
\]
Since $D\leq1$ and $p/2>(p-2)/4$, these prove the first estimate in \eqref{shell-interpolation-bounds}. For the second, use
\[
 \int_{B_R(y)}|v|^2
 \leq2|B_R|\|v_{\mathcal A}\|_\infty^2+2\|v_{\mathcal A^c}\|_2^2
 \leq C_RD^{1/2},
\]
uniformly in the center $y$. The same decomposition gives, for every $2<q\leq2_s^*$ when $N>2s$, and every finite $q>2$ when $N\leq2s$,
\begin{equation}\label{green-shell-q}
 \|v\|_q^q\leq C_q\bigl(D(v)^{(q-2)/4}+D(v)^{q/2}\bigr).
\end{equation}
Indeed, interpolate the $L^2$ and $L^\infty$ bounds for $v_{\mathcal A}$, and apply the continuous Sobolev embedding to $v_{\mathcal A^c}$. The latter embedding includes $q=2_s^*$ when that endpoint is finite.
\end{proof}

\begin{theorem}[Small-mass energy, multipliers, and delocalization]\label{smallmass-spectral}
Let $N\geq1$, $0<s<1$, and $2<p\leq2+4s/N$. For each $c>0$, let $u_c$ be any global minimizer of $\mathcal M_p(c)$, and let $\lambda_c$ be its multiplier in \eqref{variational-equation}. Put
\[
 \begin{gathered}
 b_s=\frac1{es},\qquad r_s=e^{-1/(2s)},\qquad v_c=\frac{u_c}{c},\\
\widetilde a=\frac{(N+1)(p-2)}4<2,\qquad \widetilde\nu=2+\frac{2(p-2)}{2-\widetilde a},\\
 \rho_p=\frac{4(p-2)}{6-p}>0.
 \end{gathered}
\]
The assumed range gives $p<6$ and $\widetilde a\leq s(N+1)/N<2$. There are constants $c_0,\kappa,C>0$, independent of the choice of minimizer, such that for $0<c<c_0$, 
\begin{align}
 \newadd{\kappa c^{\widetilde\nu}}&\leq-b_sc^2-\mathcal M_p(c)\leq Cc^{2+\rho_p},\label{smallmass-energy-gap}\\
 \newadd{2\kappa c^{\widetilde\nu-2}}&\leq\lambda_c-2b_s\leq Cc^{\rho_p}.\label{smallmass-multiplier-gap}
\end{align}
In particular $\mathcal M_p(c)/c^2\to-b_s$ and $\lambda_c\to2b_s$ as $c\downarrow0$. The family $(v_c)_{0<c<c_0}$ is bounded in $\HH$ and satisfies
\begin{equation}\label{smallmass-fourier-gap}
 0\leq\int_{\R^N}(\ell(\xi)+b_s)|\widehat v_c(\xi)|^2\,d\xi
       \leq Cc^{\rho_p}.
\end{equation}
Consequently, for every $\epsilon>0$,
\begin{equation}\label{smallmass-shell}
 \int_{\{\xi:\,||\xi|-r_s|\geq\epsilon\}}
          |\widehat v_c(\xi)|^2\,d\xi\longrightarrow0.
\end{equation}
In physical space these normalized profiles spread out: for every fixed $R>0$ there is $C_R$ such that
\begin{equation}\label{smallmass-spreading}
 \sup_{y\in\R^N}\int_{B_R(y)}|v_c|^2\,dx
       \leq C_R c^{\rho_p/2}\longrightarrow0.
\end{equation}
Moreover, for $2<q\leq2_s^*$ when $N>2s$, and every finite $q>2$ when $N\leq2s$,
\begin{equation}\label{green-quantitative-q}
 \|v_c\|_q^q\leq C_qc^{\rho_p(q-2)/4}\longrightarrow0.
\end{equation}
In particular, strong convergence includes the finite Sobolev endpoint. In dimension one the upper and lower powers in the two gap estimates coincide, as recorded in Corollary~\ref{one-dimensional-gaps}. No matching exponent in higher dimensions is asserted.
\end{theorem}
\begin{proof}
\emph{Uniform bounds after amplitude normalization.}
Write $K_c=K(v_c)$ and $\delta_c=c^{p-2}$. By \eqref{review-GN} and the inclusion into $H^s$,
\[
 P(v_c)\leq C(1+K_c)^r\leq C(1+K_c),\qquad
 r=\frac{N(p-2)}{4s}\leq1.
\]
Using \eqref{review-Acoercive} and
$\mathcal M_p(c)/c^2=A(v_c)-\delta_c P(v_c)/p<-b_s$,
we obtain
\[
 -b_s>\frac{\mathcal M_p(c)}{c^2}
 \geq\left(\frac12-C\delta_c\right)K_c-D-C\delta_c.
\]
For small $c$ this bounds $K_c$, hence $\|v_c\|_{\HH}$ and $P(v_c)$ uniformly. The spectral lower bound now gives
\[
 0\leq D_c:=A(v_c)+b_s
 <\frac{\delta_c}{p}P(v_c)\leq C\delta_c.
\]
This gives preliminary bounds $D_c\leq Cc^{p-2}$ and $-b_sc^2-\mathcal M_p(c)\leq Cc^p$, which will be sharpened below. Testing the equation, or using \eqref{equ14}, gives
\[
 \lambda_c-2b_s=\delta_cP(v_c)-2D_c
 =\left(1-\frac2p\right)\delta_cP(v_c)
       +\frac{2(-b_sc^2-\mathcal M_p(c))}{c^2}>0.
\]
It is at most $C\delta_c$, proving both limits and a preliminary upper multiplier bound.

\emph{An anisotropic quantitative trial bound.
Fix $\xi_0=r_se_1$ and write $\eta=(\eta_1,\eta')$. Near $\xi_0$, the identity
\[
 |\xi_0+\eta|-r_s
 =\frac{2r_s\eta_1+\eta_1^2+|\eta'|^2}{|\xi_0+\eta|+r_s}
\]
and the quadratic minimum of the radial symbol imply
\[
 0\leq\ell(\xi_0+\eta)+b_s
 \leq C\bigl(\eta_1^2+|\eta'|^4\bigr).
\]
Choose a real Schwartz function $\phi$ of unit $L^2$-norm whose Fourier transform is even, smooth, and compactly supported. Set
\[
 u_R(x)=\sqrt2\,cR^{-(N+1)/4}
 \phi\!\left(\frac{x_1}{R},\frac{x'}{\sqrt R}\right)\cos(r_sx_1).
\]
When $N=1$, the primed variables are omitted. For large $R$, the two Fourier packets are disjoint and $\|u_R\|_2=c$ exactly. Their normal and tangential widths are $O(R^{-1})$ and $O(R^{-1/2})$, respectively. Therefore, for constants $L,D_0,R_0>0$ independent of $c$,
\[
 A(u_R)\leq-b_sc^2+Lc^2R^{-2},\qquad
 P(u_R)\geq D_0c^pR^{-\widetilde a}\quad(R\geq R_0).
\]
For the nonlinear estimate, change variables and average the factor $|\cos(r_sRy_1)|^p$ as in Lemma~\ref{review-trial}. Choose $L_0$ large enough that
$LL_0^{-2}\leq D_0L_0^{-\widetilde a}/(2p)$ and put
$R=L_0c^{-(p-2)/(2-\widetilde a)}$. For small $c$ this gives
\[
 \mathcal M_p(c)\leq-b_sc^2-\frac{D_0}{2p}L_0^{-\widetilde a}c^{\widetilde\nu}.
\]
This proves the lower energy bound, and the preceding exact multiplier formula gives the lower multiplier bound. If $N>1$, then $\widetilde a<N(p-2)/2$, so the new exponent is strictly smaller than the isotropic exponent $2+2(p-2)/(2-N(p-2)/2)$ and yields a stronger lower bound. In dimension one the two exponents agree.}

\emph{A sharper bound by interpolation on the Fourier shell.}
The preliminary bound makes $D_c\leq1$ for small $c$. Also $D_c>0$, since $\ell+b_s$ vanishes only on a sphere of Lebesgue measure zero and $\|v_c\|_2=1$. Lemma~\ref{shell-interpolation}, with $\theta=(p-2)/4<1$, yields
\[
 D_c<\frac{\delta_c}{p}P(v_c)\leq C\delta_c D_c^\theta.
\]
Dividing by $D_c^\theta$ gives
\[
 D_c\leq C\delta_c^{1/(1-\theta)}=Cc^{\rho_p},\qquad
 \delta_cP(v_c)\leq C\delta_cD_c^\theta\leq Cc^{\rho_p}.
\]
The exact energy and multiplier formulas therefore give the upper bounds in \eqref{smallmass-energy-gap}--\eqref{smallmass-multiplier-gap}, together with \eqref{smallmass-fourier-gap}.

\emph{Fourier concentration and spatial spreading.}
The continuous function $\ell+b_s$ vanishes precisely on $|\xi|=r_s$ and tends to infinity at infinity. Its minimum away from any fixed neighborhood of that sphere is positive. Equation~\eqref{smallmass-fourier-gap} therefore proves \eqref{smallmass-shell}. The second estimate of Lemma~\ref{shell-interpolation} gives
\[
 \sup_y\int_{B_R(y)}|v_c|^2\leq C_RD_c^{1/2}\leq C_Rc^{\rho_p/2},
\]
which proves \eqref{smallmass-spreading}.

\emph{Quantitative vanishing, including the finite endpoint.
Apply \eqref{green-shell-q} with $D(v_c)=D_c\leq Cc^{\rho_p}$. For $0<c<1$, the term $D_c^{q/2}$ is bounded by a constant times $c^{\rho_p(q-2)/4}$, since $q/2>(q-2)/4$. This proves \eqref{green-quantitative-q}. In particular $P(v_c)\to0$. All constants are uniform in the choice of minimizer, so the conclusions hold for arbitrary minimizing choices.}
\end{proof}

\begin{corollary}[Small-mass limits in the physical normalization]\label{smallmass-physical}
In the same subcritical and critical range, let $w_c$ minimize the physical energy $I$ on $\|w_c\|_2=c$, and let $\mu_c$ be its multiplier in \eqref{equ3}. Then
\[
 \frac{I(w_c)}{c^2}\longrightarrow-\frac1{2es},\qquad
 \mu_c\longrightarrow\frac1{es}\quad\text{from above}.
\]
The Fourier-shell concentration, spatial spreading, and $L^q$ vanishing of Theorem~\ref{smallmass-spectral} hold for $w_c/c$ as well.
\end{corollary}
\begin{proof}
Use $w_c=k_pu_{c/k_p}$, $k_p=2^{-1/(p-2)}$, and
$I(k_pu)=k_p^2\E(u)/2$, $\mu_c=\lambda_{c/k_p}/2$.
Then $w_c/c=u_{c/k_p}/(c/k_p)$, so the normalized-profile conclusions are unchanged.
\end{proof}

\begin{corollary}[Matching gap powers in dimension one]\label{one-dimensional-gaps}
Let $N=1$, $0<s<1$, and $2<p\leq2+4s$. Put $\rho_p=4(p-2)/(6-p)$. There are positive constants $k_1,k_2$ such that, for every sufficiently small $c$ and every minimizing choice,
\[
 \begin{aligned}
 k_1c^{2+\rho_p}&\leq-\frac{c^2}{es}-\mathcal M_p(c)
                         \leq k_2c^{2+\rho_p},\\
 k_1c^{\rho_p}&\leq\lambda_c-\frac2{es}\leq k_2c^{\rho_p}.
 \end{aligned}
\]
The same powers hold for the physical energy gap
$-c^2/(2es)-\inf_{\|w\|_2=c}I(w)$ and the physical multiplier gap $\mu_c-1/(es)$, with possibly different positive constants. These comparisons determine the decay exponents; they do not identify a leading coefficient or a limiting rescaled envelope.
\end{corollary}
\begin{proof}
In dimension one, $\widetilde a=(p-2)/2$, and hence $\widetilde\nu=2+\rho_p$. Thus the lower and upper powers in Theorem~\ref{smallmass-spectral} coincide. The amplitude conversion in Corollary~\ref{smallmass-physical} multiplies the gaps by fixed positive constants and rescales the mass by the fixed factor $k_p$, preserving these exponents.
\end{proof}

\begin{corollary}[Quantitative convergence of the fractional kinetic moment]\label{green-kinetic-moment}
Under the assumptions of Theorem~\ref{smallmass-spectral}, put $a=N(p-2)/2$ and
\[
 D_c=\int_{\R^N}(\ell(\xi)+b_s)|\widehat v_c(\xi)|^2\,d\xi.
\]
Every global minimizer satisfies $\mathcal Q(u_c)=0$. Consequently,
\begin{equation}\label{green-moment-identity}
 \int_{\R^N}|\xi|^{2s}|\widehat v_c(\xi)|^2\,d\xi
 =\frac1e+\frac{a}{2p}c^{p-2}\|v_c\|_p^p-sD_c
 =\frac1e+O(c^{\rho_p})
 \quad(c\downarrow0),
\end{equation}
uniformly in the choice of minimizer. The same convergence rate holds for the mass-normalized physical minimizers.
\end{corollary}
\begin{proof}
Since $(u_c)_t\in\mathcal S_c$ for every $t>0$, global minimality implies that the scalar function $t\mapsto\mathcal E_p((u_c)_t)$ is minimized at $t=1$. Its explicit fiber formula is differentiable, so its derivative vanishes there, giving $\mathcal Q(u_c)=0$. No differentiability of the dilation curve in $\HH$ is required. Divide this identity by $c^2$, use $A(v_c)=-1/(es)+D_c$ and
\[
 C(N,s)\mathcal J_s(v_c,v_c)
 =2\int |\xi|^{2s}|\widehat v_c|^2,
\]
and solve for the displayed moment. Theorem~\ref{smallmass-spectral} gives $D_c=O(c^{\rho_p})$. Lemma~\ref{shell-interpolation} also gives
\[
 c^{p-2}\|v_c\|_p^p
 \leq Cc^{p-2+\rho_p(p-2)/4}=Cc^{\rho_p}.
\]
This proves the rate. The fixed amplitude and norm conversion in Corollary~\ref{smallmass-physical} preserves the exponent.
\end{proof}

\begin{corollary}[Exact energy--moment identity at the critical power]\label{critical-energy-moment}
Let $p=2+4s/N$, and let $u_c$ be any global minimizer of $\mathcal M_p(c)$. Then
\begin{equation}\label{critical-moment-exact}
 \int_{\R^N}|\xi|^{2s}\left|\widehat{u_c/c}(\xi)\right|^2\,d\xi
 =-\frac{s\mathcal M_p(c)}{c^2}>\frac1e.
\end{equation}
In particular, with the exponents of Theorem~\ref{smallmass-spectral},
\[
 s\kappa c^{\widetilde\nu-2}
 \leq\int |\xi|^{2s}\left|\widehat{u_c/c}\right|^2-\frac1e
 \leq Cc^{\rho_p}
 \qquad(0<c<c_0).
\]
For a physical minimizer $w_c$ of $I$ at norm $c$, the corresponding exact identity is
\[
 \int |\xi|^{2s}\left|\widehat{w_c/c}\right|^2
 =-\frac{2sI(w_c)}{c^2}>\frac1e.
\]
Thus the kinetic moment approaches $1/e$ from above. In dimension one its gap has the same matching power $c^{\rho_p}$ as the normalized energy gap.
\end{corollary}
\begin{proof}
At $a=N(p-2)/2=2s$, the identity $\mathcal Q(u_c)=0$ gives
\[
 2sA(u_c)+2q(u_c)-\frac{2s}{p}P(u_c)=0,
 \qquad q(u_c)=\int |\xi|^{2s}|\widehat u_c|^2.
\]
Hence $q(u_c)=-s\mathcal E_p(u_c)$. Divide by $c^2$ and use the strict spectral comparison \eqref{review-strict-bottom}. Subtracting $1/e$ and applying \eqref{smallmass-energy-gap} gives the bounds. For the physical normalization, either use the amplitude conversion or differentiate $I((w_c)_t)$ at its global minimum $t=1$ to obtain $q(w_c)=-2sI(w_c)$. Corollary~\ref{one-dimensional-gaps} gives the matching one-dimensional power.
\end{proof}

\section{Normalized Solutions in the Mass-Critical Case}

\begin{lemma}\label{lemma-critical-coercivity}
At $p=2+4s/N$, the functional $\E$ is bounded below and coercive on $\Sc$ for every $c>0$, and $\mathcal M_p(c)<-c^2/(es)$.
\end{lemma}
\begin{proof}
Use \eqref{review-critical-coercive} and Lemma~\ref{review-trial}.
\end{proof}

\begin{proof}[Proof of Theorem~\ref{theorem3}]
The proof of Theorem~\ref{theorem1}, including strict subadditivity and exact mass recovery, applies at $p=2+4s/N$. The exponent remains strictly below $2_s^*$ whenever that endpoint is finite. This proves existence, compactness, and the positive multiplier assertion at every mass.
\end{proof}

\begin{remark}[Effect of the logarithmic correction]
Here the term mass-critical refers to the fractional power $2+4s/N$. The logarithmic equation is not scale invariant at that power. The estimate $S\leq\delta K+D_\delta\|u\|_2^2$ permits absorption of the critical nonlinearity for every fixed mass, which explains the absence of a finite critical-mass obstruction to global minimization. This conclusion is specific to the strengthened logarithmic kinetic energy.
\end{remark}

\section{Normalized Solutions in the Mass-Supercritical Case}
\subsection{Notation and exact fiber geometry}
Throughout, $0<s<1$, $2+\frac{4s}{N}<p<2_s^*$, and
\[
 a=\frac{N(p-2)}2,\qquad d=a-2s>0,\qquad C=C(N,s)>0,\qquad B=B(N,s).
\]
For $u\ne0$ put
\[
 A=A(u)=\mathcal E_{s+\log}(u,u),\quad
 b=b(u)=C\Js(u,u)>0,\quad P=P(u)=\|u\|_p^p>0.
\]
We use exactly the energy normalization in the definition above:
\begin{equation}\label{fiber}
 f_u(t):=\E(u_t)=t^{2s}(A+b\ln t)-\frac{P}{p}t^a,
 \qquad u_t(x)=t^{N/2}u(tx).
\end{equation}
Thus
\begin{equation}\label{G}
 G_u(t):=t^{-2s}\Q(u_t)
 =2sA+b+2sb\ln t-\frac ap Pt^d,
 \qquad f_u'(t)=t^{2s-1}G_u(t).
\end{equation}
The mass convention is $\|u\|_2=c$, so $\|u\|_2^2=c^2$.
The set $\Bc$ in the definition above is characterized by
\begin{equation}\label{branch}
 2sb<\frac{ad}{p}P.
\end{equation}
In particular, $\Pc$ consists of the zeros of $\Q$ that satisfy this strict inequality. It is the branch of strict local fiber maxima, rather than all the zeros of $\Q$.

\begin{fibercriterion}
Define
\begin{equation}\label{tstar}
 t_* =\left(\frac{2sbp}{adP}\right)^{1/d}.
\end{equation}
Then $G_u$ is strictly increasing on $(0,t_*)$ and strictly decreasing on $(t_*,\infty)$, and tends to $-\infty$ at both endpoints. A scale $t>0$ with $u_t\in\Pc$ exists if and only if $G_u(t_*)>0$. In that case there are exactly two zeros
\[
 0<t_-<t_*<t_+,
\]
and only $u_{t_+}$ belongs to $\Pc$. If $G_u(t_*)=0$, the unique zero is degenerate and does not belong to $\Pc$; if $G_u(t_*)<0$, there is no zero.
\end{fibercriterion}
\begin{proof}
Direct differentiation gives
\[
 G_u'(t)=\frac1t\left(2sb-\frac{ad}{p}Pt^d\right).
\]
This proves the monotonicity and the root count. At a zero, the strict branch condition for $u_t$ is precisely $t>t_*$. Moreover,
\[
 f_u''(t)=t^{2s-1}G_u'(t)\quad\hbox{when }G_u(t)=0,
\]
Thus $f_u$ has a strict local minimum at the smaller zero and a strict local maximum at the larger zero.
\end{proof}

\begin{lemma}\label{lemma3}
Let $0<c<c_*$, with $c_*$ as in \eqref{newconstants}. Any minimizer of $\E$ on $\Pc$ is a constrained critical point of $\E$ on $\Sc$.
\end{lemma}
\begin{proof}
The uniform fiber barrier and the smooth reduced functional $\Psi_c$ are established below. If $u\in\Pc$ minimizes $\E$, then $t_+(u)=1$ and $u$ minimizes $\Psi_c$ on $\Sc$. The scalar chain rule, together with $\Q(u)=0$, yields $D\Psi_c(u)[h]=D\E(u)[h]$ for every tangent vector $h$ to $\Sc$. Hence $D(\E|_{\Sc})(u)=0$. This argument does not presume that the dilation orbit is differentiable in the energy norm.
\end{proof}

\subsection{The fiber remainder and minimax characterization}
\begin{lemma}[Fiber-remainder identity]\label{lemma4}
For every $u\in\HH$ and $t>0$,
\begin{equation}\label{identity}
 \E(u)=\E(u_t)+\frac{1-t^{2s}}{2s}\Q(u)+\mathcal H(t,u),
\end{equation}
where
\begin{equation}\label{H}
 \mathcal H(t,u)=b\left(\frac{t^{2s}-1}{2s}-t^{2s}\ln t\right)
 +\frac Pp\left(t^a-1+\frac a{2s}(1-t^{2s})\right).
\end{equation}
For each $u\in\Bc$ there is a number $\delta(u)\in(0,1)$ such that
\[
 \mathcal H(1,u)=0,\qquad
 \mathcal H(t,u)>0\quad
 \text{for }t\in(1-\delta(u),1)\cup(1,\infty).
\]
No uniform choice of $\delta$ over all of $\Bc$ is asserted.
\end{lemma}
\begin{proof}
The identity follows by substituting \eqref{fiber} and $\Q(u)=2sA+b-aP/p$. Differentiation yields
\begin{equation}\label{Hprime}
 \partial_t\mathcal H(t,u)=t^{2s-1}F_u(t),\qquad
 F_u(t)=-2sb\ln t+\frac ap P(t^d-1).
\end{equation}
This expression has $F_u(1)=0$. In particular,
\[
 F_u(1)=0,\qquad F_u'(1)=-2sb+\frac{ad}{p}P>0.
\]
Continuity implies $F_u(t)<0$ immediately to the left of $1$, which gives the claimed local positivity of $\mathcal H$ there. For $t>1$, the inequality $t^d-1>d\ln t$ gives
\[
 F_u(t)>\left(-2sb+\frac{ad}{p}P\right)\ln t>0.
\]
Integrating from $1$ proves positivity on $(1,\infty)$.
\end{proof}

The dependence on $u$ matters. A uniform left neighborhood cannot be obtained just from a strict inequality with no uniform margin. At the limiting ratio $adP/p=2sb$, writing $t=e^x$ gives
\[
 F_u(e^x)=\frac{2sb}{d}\bigl(e^{dx}-1-dx\bigr)>0\quad(x<0).
\]
Consequently $\mathcal H(t,u)<0$ for $t<1$ at that limiting ratio. The strict-branch ratio can approach this value: for a fixed smooth nonzero profile, increasing its mass-preserving dilation makes $P/b$ range over $(0,\infty)$.

\begin{corollary}[One-sided minimax characterization]\label{corollary1}
If $\Pc\ne\varnothing$, then
\begin{equation}\label{minimaxeq}
 \gamma_p(c)=\inf_{u\in\widetilde{\mathbf S}_c}
       \max_{t\geq1}\E(u_t),\qquad
 \widetilde{\mathbf S}_c=\{u\in\Bc:\Q(u)\geq0\}.
\end{equation}
For each $u\in\Pc$, the fiber $t\mapsto\E(u_t)$ also has a strict local maximum at $t=1$.
\end{corollary}
\begin{proof}
If $u\in\widetilde{\mathbf S}_c$, then $t_*<1$, $G_u(1)\geq0$, and $G_u$ is strictly decreasing on $[1,\infty)$. There is a unique $t_+\geq1$ with $G_u(t_+)=0$. The fiber energy increases up to $t_+$ on this interval and decreases thereafter. Hence
\[
 \max_{t\geq1}\E(u_t)=\E(u_{t_+})\geq\gamma_p(c).
\]
Conversely, for $u\in\Pc$ we have $t_+=1$, and the maximum is $\E(u)$. Taking infima proves \eqref{minimaxeq}. The local assertion follows from Lemma~\ref{lemma4}, or from $f_u''(1)<0$.
\end{proof}
This formulation avoids a common left endpoint for all fibers. The global maximum on $t>0$ need not be achieved at $t_+$, since $\E(u_t)\to0$ as $t\downarrow0$, while $\E(u_{t_+})$ may be negative.

\subsection*{Projection when \texorpdfstring{$\Q(u)\leq0$}{the \Pohozaev{} functional is nonpositive}}
\begin{lemma}[Projection onto the upper branch]\label{lemma5}
Suppose $u\in\Bc$. If $\Q(u)\geq0$, there is a unique $t_+\geq1$ such that $u_{t_+}\in\Pc$, with $t_+=1$ exactly when $\Q(u)=0$.

If $\Q(u)\leq0$, such a projection exists if and only if
\begin{equation}\label{barrier}
 G_u(t_*)=2sA+b+\frac{2sb}{d}
 \left[\ln\left(\frac{2sbp}{adP}\right)-1\right]>0.
\end{equation}
When \eqref{barrier} holds, the projection is unique on the upper branch and
\[
 t_*<t_+\leq1,
\]
with $t_+<1$ exactly when $\Q(u)<0$.
\end{lemma}
\begin{proof}
These assertions follow from the fiber-root criterion, $t_*<1$, and the strict decrease of $G_u$ on $(t_*,\infty)$. Formula \eqref{barrier} follows by substituting \eqref{tstar} in \eqref{G}.
\end{proof}
\noindent Thus the exponent range alone, even together with $u\in\Bc$ and $\Q(u)<0$, does not justify projection. One must check the positive fiber barrier \eqref{barrier}. In particular, replacing $\Q\geq0$ by $\Q\leq0$ without the barrier condition would be incorrect.

For a concrete obstruction, fix $w\in C_c^\infty(\R^N)\setminus\{0\}$ and put $u=qw$, $q>0$. With $A_0=A(w)$, $b_0=b(w)$, and $P_0=P(w)$,
\[
 \frac{G_{qw}(t_*(qw))}{q^2}
 =2sA_0+b_0+\frac{2sb_0}{d}
 \left[\ln\left(\frac{2sb_0p}{adP_0}\right)-(p-2)\ln q-1\right].
\]
For large $q$, this quantity is negative, $qw\in\mathbf S_{q\|w\|_2}$, and $\Q(qw)<0$. Nevertheless, its entire fiber has no Poho\u{z}aev point.

\subsection{Positivity and the small-mass restriction}
On $\Pc$, eliminating $A$ or $P$ gives the two useful exact identities
\begin{equation}\label{energyP}
 \E(u)=\frac{d}{2sp}P-\frac{b}{2s}
       =\frac da A-\frac ba.
\end{equation}
Therefore
\begin{equation}\label{positivecriterion}
 \E(u)>0\quad\Longleftrightarrow\quad P>\frac{pb}{d}.
\end{equation}
The branch condition only requires $P>2spb/(ad)$. Since $a>2s$, it is strictly weaker than \eqref{positivecriterion}.

\begin{proposition}[Negative energy on the strict upper branch]\label{negative}
For every admissible choice of $N$, $s$, and $p$, there exist $c>0$ and a function $u\in\Pc$ such that $\E(u)<0$. In particular, whenever $\gamma_p(c)$ is defined as the infimum over this set, $\gamma_p(c)\leq\E(u)<0$ at that mass.
\end{proposition}
\begin{proof}
Fix a smooth nonzero $w$ and use $A_0,b_0,P_0$ as above. Choose
\begin{equation}\label{tangent}
 T=\exp\left(\frac1d-\frac{2sA_0+b_0}{2sb_0}\right),\qquad
 q_*^{p-2}=\frac{2sb_0p}{adP_0T^d}.
\end{equation}
For $q=q_*$, the normalized Poho\u{z}aev function $G_{q_*w}$ has its unique maximum at $T$, where it equals zero. Indeed,
\[
 G_{q_*w}(T)=0,\qquad G'_{q_*w}(T)=0.
\]
At $v_*=(q_*w)_T$, the equality $a d P(v_*)/p=2sb(v_*)$ and \eqref{energyP} give
\begin{equation}\label{negativevalue}
 \E(v_*)=-\frac{d}{2sa}\,b(v_*)<0.
\end{equation}
Although $v_*$ itself is degenerate and excluded from $\Pc$, it is a limit of strict upper-branch points. To see this, for $0<q<q_*$ the maximum of $G_{qw}$ is
\[
 G_{qw}(t_*(qw))=q^2\frac{2sb_0(p-2)}d\ln\frac{q_*}{q}>0.
\]
There are two roots, and the upper one $t_+(q)$ tends to $T$ as $q\uparrow q_*$. For completeness, the roots stay in a compact subinterval of $(0,\infty)$ when $q$ is close to $q_*$, by the uniform endpoint limits of $G_{qw}/q^2$; every limit root must be the unique root $T$ at $q_*$. Consequently,
\[
 v_q=(qw)_{t_+(q)}\in\mathcal P_{q\|w\|_2},\qquad
 \E(v_q)\longrightarrow\E(v_*)<0.
\]
For $q<q_*$ sufficiently close to $q_*$, take $c=q\|w\|_2$ and $u=v_q$.
\end{proof}

\subsection*{Lower bounds at small mass}
We next establish positivity for sufficiently small mass, provided $\Pc$ is nonempty. The subsequent compactness argument supplies attainment.

Write $S(u)=\Js(u,u)$, $K(u)=\Jp(u,u)$, and $\newadd{\sigma_{N-1}=|\mathbb S^{N-1}|}$. The elementary inequality
\begin{equation}\label{infinitesimal}
 S(u)\leq\varepsilon K(u)+D_\varepsilon\|u\|_2^2
 \qquad(\varepsilon>0)
\end{equation}
holds with a finite $D_\varepsilon$. Indeed, for $0<r<1$, split the difference integral at $|x-y|=r$. The near part is at most $K(u)/(C\ln(1/r))$, while the far part is at most
\[
 \frac{2\newadd{\sigma_{N-1}}}s r^{-2s}\|u\|_2^2.
\]
Choosing $r$ sufficiently small proves \eqref{infinitesimal}. Similarly,
\begin{equation}\label{tail}
 \Jm(u,u)\leq M\|u\|_2^2,\qquad
 M=4C\newadd{\sigma_{N-1}}\int_1^\infty(\ln r)r^{-1-2s}\,dr
 =\frac{C\newadd{\sigma_{N-1}}}{s^2}.
\end{equation}
Since $A=K-\Jm+(BC/2)S$, equations \eqref{infinitesimal}--\eqref{tail} imply that, for every $L>0$, some finite $D_L>0$ satisfies
\begin{equation}\label{Alower}
 A(u)\geq L S(u)-D_L\|u\|_2^2.
\end{equation}
For example, take $\varepsilon>0$ with $\varepsilon^{-1}+BC/2\geq L$ and then $D_L=D_\varepsilon/\varepsilon+M$.

\begin{lemma}[Lower bound and positivity at small mass]\label{lemma7}
Set
\[
 r_0=\frac{a}{2s}>1,\qquad \beta=p-\frac as>0.
\]
Let $G_{N,s,p}$ be a constant in the usual fractional interpolation inequality
\begin{equation}\label{usualGN}
 \|u\|_p^p\leq G_{N,s,p} S(u)^{r_0}\|u\|_2^{\beta}.
\end{equation}
Then, for all $u\in\Bc$,
\begin{equation}\label{rho}
 S(u)>\rho_c:=R c^{-\beta/(r_0-1)},\qquad
 R=\left(\frac{2sCp}{adG_{N,s,p}}\right)^{1/(r_0-1)}.
\end{equation}
Choose $D=D_{2C/d}>0$ in \eqref{Alower}. If
\begin{equation}\label{masssmall}
 0<c\leq c_0:=
 \left(\frac{CR}{2dD}\right)^{1/(2+\beta/(r_0-1))},
\end{equation}
then for every $u\in\Pc$,
\begin{equation}\label{coercive}
 \E(u)\geq\frac{C}{2a}S(u)\geq\frac{C}{2a}\rho_c>0.
\end{equation}
Consequently $\gamma_p(c)>0$ whenever $0<c\leq c_0$ and $\Pc\ne\varnothing$.
\end{lemma}
\begin{proof}
Combining \eqref{branch} and \eqref{usualGN} gives
\[
 2sC S(u)<\frac{ad}{p}G_{N,s,p}S(u)^{r_0}c^\beta,
\]
which proves \eqref{rho}. Next use \eqref{energyP} and \eqref{Alower} with $L=2C/d$:
\[
 \E(u)\geq\frac Ca S(u)-\frac{dD}{a}c^2.
\]
The mass condition \eqref{masssmall} is equivalent to $dDc^2\leq C\rho_c/2$. Since $S(u)\geq\rho_c$, this proves \eqref{coercive}.
\end{proof}

\subsection{Proof of the small-mass existence theorem}\label{replacement}
We retain the normalization $\E(u)=A(u)-\|u\|_p^p/p$ used above. Let $\mathscr L=(-\Delta)^{s+\log}$ with Fourier symbol
\[
 \ell(\xi)=|\xi|^{2s}\ln|\xi|^2.
\]
The unitary Fourier transform is used, with $(-\Delta)^s$ having symbol $|\xi|^{2s}$. In the present notation,
\begin{equation}\label{fourierexact}
 A(u)=\int_{\R^N}\ell(\xi)|\widehat u(\xi)|^2\,d\xi,
 \qquad b(u)=2\int_{\R^N}|\xi|^{2s}|\widehat u(\xi)|^2\,d\xi.
\end{equation}
These identities also follow by differentiating the fractional quadratic form in $s$ and then using density. All spaces and variations below may be taken over the real numbers.

\subsubsection*{A quantitative mass range and the theorem}
Continue to write
\begin{equation}\label{newconstants0}
 r_0=\frac a{2s}>1,\quad \beta=p-\frac as>0,\quad
 \sigma=\frac{\beta}{r_0-1},\quad \rho_m=R m^{-\sigma},
\end{equation}
with $R$ as in \eqref{rho}. Put
\begin{equation}\label{newconstants}
 \begin{split}
 K_b&=G_{N,s,p}C^{-r_0},\\
 c_B&=\left(\frac{p}{\max\{a,4s\}K_b(2e)^{r_0-1}}\right)^{1/(p-2)},\\
 \Lambda_{\rm low}&=2\exp(2s/d-1)
       \left(\frac{2s\beta}{ad}+\frac1s\right),\\
 \Lambda&=1+\max\left\{\Lambda_{\rm low},\frac{2}{es}\right\},\qquad
 c_\lambda=\left(\frac{CR}{2s\Lambda}\right)^{1/(\sigma+2)},\\
 c_*&=\frac12\min\{c_0,c_B,c_\lambda\}>0,
 \end{split}
\end{equation}
where $c_0$ is the positivity threshold in \eqref{masssmall}. These constants give sufficient mass restrictions; we do not claim that they are optimal.

The following arguments prove Theorem~\ref{theorem4}.

\subsubsection*{The uniform fiber barrier}
For $u\in\mathcal S_m$, Jensen's inequality applied to $y\ln y$ and the probability measure $|\widehat u|^2/m^2$ gives
\begin{equation}\label{jensen}
 A(u)\geq\frac{b(u)}{2s}\ln\frac{b(u)}{2m^2}.
\end{equation}
Indeed, with $y=|\xi|^{2s}$ we have $\ell(\xi)=y\ln y/s$ and $\int y|\widehat u|^2=b(u)/2$. The negative part of $y\ln y$ is bounded, so all terms in this use of Jensen are legitimate.

On each mass-preserving fiber, choose $t_b$ such that $b(u_{t_b})=2em^2$. Combining \eqref{jensen} with
\begin{equation}\label{GNb}
 P(v)\leq K_b b(v)^{r_0}m^\beta\qquad(v\in\mathcal S_m)
\end{equation}
and $2r_0+\beta=p$ yields, for $m<c_B$,
\begin{equation}\label{barrieruniform}
 \begin{split}
 \Q(u_{t_b})&\geq 2b(u_{t_b})-\frac ap P(u_{t_b})
                       \geq b(u_{t_b})>0,\\
 \E(u_{t_b})&\geq\frac{b(u_{t_b})}{2s}-\frac1pP(u_{t_b})
                       \geq\frac{b(u_{t_b})}{4s}>0.
 \end{split}
\end{equation}
The fiber-root criterion now gives two distinct roots. Since $t_b$ lies between them, the upper fiber maximum is positive and is the global maximum over $t>0$: the other endpoint limit is $0$ and the large-$t$ limit is $-\infty$. In particular,
\begin{equation}\label{globalreduced}
 \gamma_p(m)=\inf_{u\in\mathcal S_m}\Psi_m(u),\qquad
 \Psi_m(u):=\E(u_{t_+(u)})=\max_{t>0}\E(u_t).
\end{equation}
The positive fiber barrier justifies this all-fiber formulation in the stated mass range.

\subsubsection*{Smooth reduction and strict mass comparison}
The quantities $A(u)$, $b(u)$, and $P(u)$ are continuously differentiable on $\HH$. Formula \eqref{G} and the strict upper-root condition imply, by the implicit function theorem, that $u\mapsto t_+(u)$ is $C^1$ on the mass sphere. The reduced functional $\Psi_m$ is $C^1$ as well. This conclusion uses its explicit expression in the three scalar quantities $A,b,P$; it does not require differentiability of $t\mapsto u_t$ in $\HH$.

If $t=t_+(u)$ and $h$ is tangent to the mass sphere, the cancellation of the derivative in $t$ gives
\begin{equation}\label{reducedderivative}
 D\Psi_m(u)[h]=D\E(u_t)[h_t].
\end{equation}
For a fixed profile $w$ with $\|w\|_2=1$, set
\[
 v_m=(mw)_{t_+(mw)},\qquad M_w(m)=\E(v_m),\qquad 0<m<c_*.
\]
Differentiating only the explicit scalar expression for the energy and using $\Q(v_m)=0$ gives
\begin{equation}\label{massderivative}
 M_w'(m)=\frac1m\bigl(2A(v_m)-P(v_m)\bigr)
 =-\frac1m\left(\frac\beta p P(v_m)+\frac1s b(v_m)\right)
 \leq-\frac{C\rho_m}{sm}.
\end{equation}
Integrate between $m_1$ and $m_2$, start with profiles whose upper-branch energies approach $\gamma_p(m_1)$, and use \eqref{globalreduced} at $m_2$. 
This proves \eqref{strictquant} using the convention $\|v_m\|_2=m$ and the smooth upper projection.

\subsubsection*{A bounded minimizing sequence}
Let $(z_n)\subset\Pc$ be a minimizing sequence for $\gamma_p(c)$. Lemma~\ref{lemma7} gives a uniform bound on $S(z_n)$. On $\Pc$, \eqref{energyP} implies
\[
 A(z_n)=\frac{a\E(z_n)+b(z_n)}d.
\]
Thus $A(z_n)$ is bounded. Since
\[
 K(z_n)=A(z_n)+\Jm(z_n,z_n)-\frac{BC}{2}S(z_n),
\]
the tail bound \eqref{tail} proves boundedness in $X:=\HH$.

The derivative of the normalized fiber equation at an upper root, evaluated at $t=1$, satisfies the particularly useful identity
\begin{equation}\label{gap}
 D(z):=2sb(z)-\frac{ad}{p}P(z)
 =-d\,b(z)-2sa\E(z)\leq-dC\rho_c<0
 \qquad(z\in\Pc).
\end{equation}
Apply Ekeland's variational principle to the $C^1$ functional $\Psi_c$, which is bounded below on the complete metric space $\Sc$, starting with the sequence $(z_n)$. We obtain $\widetilde z_n\in\Sc$ such that
\[
 \|\widetilde z_n-z_n\|_X\to0,\qquad
 \Psi_c(\widetilde z_n)\to\gamma_p(c),\qquad
 \|D(\Psi_c|_{\Sc})(\widetilde z_n)\|\to0.
\]
The scalar upper projection is stable uniformly along this bounded minimizing sequence. Indeed, $A$, $b$, and $P$ are continuous on $X$, so their values at $\widetilde z_n$ and $z_n$ differ by $o(1)$. The coefficient triples at $z_n$ remain in a compact subset of $\mathbb R^3$. At $t=1$ their fiber functions vanish, and their derivatives are at most $-dC\rho_c$ by \eqref{gap}. Uniform continuity of the scalar derivatives near these triples gives a common neighborhood of $t=1$ on which they remain strictly negative. The implicit function theorem, or the mean value theorem on this neighborhood, therefore gives
\[
 \tau_n:=t_+(\widetilde z_n)\longrightarrow1.
\]
Set $u_n=(\widetilde z_n)_{\tau_n}\in\Pc$. The operators $h\mapsto h_t$ and their inverses are uniformly bounded on $X$ for $t$ in a compact neighborhood of $1$. For example, \eqref{Kscale} bounds $K(h_t)$ by a constant times $K(h)+S(h)$ there; use \eqref{infinitesimal} for the latter term. Consequently \eqref{reducedderivative} implies
\begin{equation}\label{PS}
 \E(u_n)\to\gamma_p(c),\quad \Q(u_n)=0,\quad
 \E'(u_n)+\lambda_n u_n\longrightarrow0\quad\text{in }X^*.
\end{equation}
Here $\lambda_n$ may be chosen by testing the residual with $u_n$:
\begin{equation}\label{lambdan}
 \lambda_n c^2=P(u_n)-2A(u_n)
 =\frac\beta p P(u_n)+\frac1s b(u_n).
\end{equation}
The sequence is bounded in $X$, so $\lambda_n$ is bounded. After extraction,
\begin{equation}\label{lambdalimit}
 \lambda_n\to\lambda\geq\frac{C\rho_c}{sc^2}
 =\frac{CR}{s}c^{-\sigma-2}>\Lambda.
\end{equation}

\subsubsection*{Poho\u{z}aev identities for weak profiles}
We record a variational argument proving the required identity for weak profiles without differentiating their dilation orbits in the energy norm.

\begin{lemma}[Variation of a single profile]\label{profilevariation}
Let $0<c<c_*$ and let $(u_n)\subset\Pc$ be any bounded minimizing sequence at $\gamma_p(c)$. If, after translations, $u_n\rightharpoonup v\ne0$ in $X$, then $\Q(v)=0$.
\end{lemma}
\begin{proof}
Work in the translated coordinates, pass to almost-everywhere convergence, and put $r_n=u_n-v$. Along a subsequence, the bounded scalar sequences $A(u_n)$, $b(u_n)$, and $P(u_n)$ converge to $A_\infty,b_\infty,P_\infty$. Their limits satisfy
\[
 2sA_\infty+b_\infty-\frac apP_\infty=0,\qquad
 A_\infty-\frac1pP_\infty=\gamma_p(c).
\]
For a fixed $t$ near $1$, replace only the profile $v$ by $v_t$ and restore the mass:
\[
 w_n(t)=\frac{c}{\|r_n+v_t\|_2}(r_n+v_t).
\]
Weak convergence and $\|v_t\|_2=\|v\|_2$ imply $\|r_n+v_t\|_2\to c$. Quadratic splitting and the Br\'ezis--Lieb lemma imply that the three scalar coefficients of $w_n(t)$ converge, respectively, to
\begin{equation}\label{modifiedcoeff}
 \begin{split}
 A_t&=A_\infty-A(v)+t^{2s}\bigl(A(v)+b(v)\ln t\bigr),\\
 b_t&=b_\infty-b(v)+t^{2s}b(v),\\
 P_t&=P_\infty-P(v)+t^aP(v).
 \end{split}
\end{equation}
At $t=1$, the upper-root derivative is bounded away from zero by \eqref{gap}. Therefore, for $t$ close to $1$, the upper projection scales of $w_n(t)$ converge to a $C^1$ scalar function $\eta(t)$, with $\eta(1)=1$. In particular,
\[
 \Psi_c(w_n(t))\longrightarrow
 F(t):=\eta(t)^{2s}\bigl(A_t+b_t\ln\eta(t)\bigr)
                  -\frac{\eta(t)^a}{p}P_t.
\]
For every such fixed $t$, \eqref{globalreduced} gives $F(t)\geq\gamma_p(c)$, while $F(1)=\gamma_p(c)$. Hence $F'(1)=0$. The term containing $\eta'(1)$ is the total \Pohozaev{} functional and vanishes. Differentiating \eqref{modifiedcoeff} now yields
\[
 0=F'(1)=2sA(v)+b(v)-\frac apP(v)=\Q(v).
\]
All differentiations here are of the explicit limiting scalar expressions; no unbounded spatial dilation generator is used.
\end{proof}

\subsubsection*{Exclusion of lower-branch profiles}
Suppose a translated weak limit $v\ne0$ of the sequence \eqref{PS} has norm $m\leq c$. Local compactness and the weak equation imply
\begin{equation}\label{profileequation}
 2\mathscr L v+\lambda v=|v|^{p-2}v.
\end{equation}
Testing with $v$ is legitimate in $X$. Lemma~\ref{profilevariation} gives $\Q(v)=0$, so
\begin{equation}\label{profilelambda}
 \lambda m^2=\frac\beta p P(v)+\frac1s b(v).
\end{equation}
If $v$ were on the lower or degenerate branch, then
\[
 \frac{ad}{p}P(v)\leq2sb(v).
\]
Combining this inequality with $\Q(v)=0$ gives
\[
 \frac{2sA(v)}{b(v)}\leq\frac{2s}{d}-1.
\]
The Jensen bound \eqref{jensen}, now at mass $m$, implies
\[
 \frac{b(v)}{m^2}\leq2\exp(2s/d-1).
\]
Substitute these estimates in \eqref{profilelambda} to obtain
\[
 \lambda\leq\Lambda_{\rm low},
\]
contradicting \eqref{lambdalimit}. Thus every nonzero profile belongs to $\mathcal P_m$. In particular, it has positive energy and satisfies $S(v)\geq\rho_m\geq\rho_c$.

\subsubsection*{Compactness and exact recovery of mass}
We give the extraction argument to distinguish genuine mass recovery from an inference based only on  nonincreasing levels.

Since $(u_n)$ is bounded in $H^s$, the usual Lions vanishing criterion applies. Vanishing would give $P(u_n)\to0$, contradicting $2sb(u_n)<adP(u_n)/p$ and $b(u_n)\geq C\rho_c$. Therefore some translations yield a nonzero weak limit $v^1$ in $X$. Local compactness gives almost-everywhere and local $L^p$ convergence. The residual after subtracting this translated profile has the usual quadratic norm splitting and \rev{Br\'ezis--Lieb} splitting of its $p$th power.

If the residual does not vanish locally in $L^2$, extract a second translated profile. Its center must diverge from the first, since the first residual converges weakly to zero in that frame. Repeat. Every extracted profile is a translated weak limit of the original sequence, solves \eqref{profileequation}, and satisfies Lemma~\ref{profilevariation}. The preceding subsection therefore shows that each consumes at least $\rho_c$ of the squared fractional seminorm. Orthogonality gives
\[
 \sum_{j=1}^k S(v^j)\leq\limsup_{n\to\infty}S(u_n)<\infty.
\]
Consequently, this procedure must stop after finitely many profiles, with a residual $r_n$ that vanishes in $L^p$. We have
\begin{equation}\label{finiteprofiles}
 u_n=\sum_{j=1}^k v^j(\,\cdot-y_n^j)+r_n,\qquad
 |y_n^i-y_n^j|\to\infty\ (i\ne j),\qquad \|r_n\|_p\to0.
\end{equation}
At each extraction, derivative splitting also holds:
\[
 |u_n|^{p-2}u_n-|u_n-v|^{p-2}(u_n-v)-|v|^{p-2}v
 \longrightarrow0\quad\text{in }L^{p/(p-1)}
\]
in that profile's frame. This follows from almost-everywhere convergence, boundedness in $L^p$, and the elementary epsilon estimates for the power map. Subtracting the finitely many profile equations from \eqref{PS}, using $\lambda_n\to\lambda$, and testing with $r_n$ gives
\begin{equation}\label{remaindercoerc}
 2A(r_n)+\lambda\|r_n\|_2^2=P(r_n)+o(1)=o(1).
\end{equation}
The Fourier symbol satisfies
\begin{equation}\label{spectralbottom}
 \min_{\xi}\ell(\xi)=-\frac1{es}.
\end{equation}
Since $\lambda>2/(es)$, the quadratic form on the left of \eqref{remaindercoerc} is positive and controls the entire $X$ norm. To check the latter claim directly, combine $A=K-\Jm+(BC/2)S$ with \eqref{infinitesimal}--\eqref{tail} to obtain $A\geq K/2-D\|u\|_2^2$, for some $D$. A sufficiently small convex combination of this bound and $2A+\lambda\|u\|_2^2\geq(\lambda-2/(es))\|u\|_2^2$ gives
\[
 2A(u)+\lambda\|u\|_2^2\geq\kappa_\lambda\|u\|_X^2
 \quad\text{for some }\kappa_\lambda>0.
\]
Thus $r_n\to0$ strongly in $X$. In particular, if $m_j=\|v^j\|_2$,
\begin{equation}\label{massenergyexact}
 \sum_{j=1}^k m_j^2=c^2,\qquad
 \gamma_p(c)=\sum_{j=1}^k\E(v^j).
\end{equation}
Each profile belongs to $\mathcal P_{m_j}$, so strict mass monotonicity and positivity give
\[
 \gamma_p(c)=\sum_{j=1}^k\E(v^j)
 \geq\sum_{j=1}^k\gamma_p(m_j)
 \geq k\gamma_p(c)>0.
\]
If $k\geq2$, this is impossible. Therefore $k=1$, and \eqref{massenergyexact} gives $m_1=c$. After one sequence of translations, $u_n\to u_c=v^1$ strongly in $X$, with $\E(u_c)=\gamma_p(c)$ and $u_c\in\Pc$. Passing to the limit in \eqref{PS} proves \eqref{equationcorrect}.

Finally, the argument applies to any initial minimizing sequence $(z_n)$ through the Ekeland construction. We have $\widetilde z_n-z_n\to0$ in $X$, $u_n=(\widetilde z_n)_{\tau_n}$, and $\tau_n\to1$. Once the projected sequence is relatively compact modulo translations, strong continuity of dilations at its limit transfers compactness back to $\widetilde z_n$ and $z_n$. This proves all the compactness assertions in Theorem~\ref{theorem4}. More explicitly, if $u_n(\cdot+y_n)\to u$ in $X$, then
\[
 \widetilde z_n(\cdot+\tau_ny_n)
 =\bigl(u_n(\cdot+y_n)\bigr)_{1/\tau_n}\longrightarrow u
 \quad\text{in }X.
\]
Here uniform boundedness of the dilation operators controls the convergent error, and strong continuity is applied to the fixed limit $u$.

\subsubsection*{The stationary equation and the energy normalization}
We now translate the conclusion into the physical energy convention. By \eqref{fourierexact}, $A(u)=\langle\mathscr L u,u\rangle$, so
\[
 D\bigl(A(u)-P(u)/p\bigr)[h]
 =2\langle\mathscr L u,h\rangle-\int |u|^{p-2}u h.
\]
Therefore \eqref{equationcorrect}, with $2\mathscr L$, is the equation directly associated with the functional $\E$. To obtain exactly
\begin{equation}\label{physicalequation}
 \mathscr L v+\mu v=|v|^{p-2}v,
\end{equation}
the consistent energy is
\[
 I(v)=\frac12A(v)-\frac1pP(v).
\]
Alternatively, apply Theorem~\ref{theorem4} and set
\begin{equation}\label{amplitudemap}
 k_p=2^{-1/(p-2)},\qquad v=k_p u_c,\qquad\mu=\lambda_c/2.
\end{equation}
This gives \eqref{physicalequation} with $\|v\|_2=k_pc$ and $\mu>1/(es)$. Moreover,
\[
 I(k_pu)=\frac{k_p^2}{2}\E(u),\qquad
 \left.\frac{d}{dt}I((k_pu)_t)\right|_{t=1}
       =\frac{k_p^2}{2}\Q(u),
\]
so the map preserves the upper-branch minimization problem and its positivity, with the rescaled mass. Consequently the original stationary equation has an upper-branch normalized minimizer for every mass
\[
 0<\|v\|_2<2^{-1/(p-2)}c_*.
\]
The positive-energy conclusion remains restricted to this upper-branch problem. Lower-branch critical points have not been excluded at the prescribed mass, so calling this solution a ground state among all normalized critical points would require a further comparison.

\subsection{Radial normalized solutions at every mass}
In this subsection we assume $N\geq2$ and allow $2+4s/N<p\leq2_s^*$, extending the exponent range used in the preceding subsections. Since $0<s<1$, we have $N>2s$ and the endpoint is finite.
We now address the existence of solutions for an arbitrary prescribed mass by a different variational construction. The low-energy level $\gamma_p(c)$ is not used. Instead, a symmetric linking construction produces high-energy constrained critical points in the radial space. All calculations in this subsection use the physical energy
\begin{equation}\label{am-energy}
 \begin{aligned}
 I(u)&=\frac12 A(u)-\frac1pP(u),\qquad
 A(u)=\langle\mathscr L u,u\rangle,\qquad P(u)=\|u\|_p^p,\\
 q(u)&=\int_{\R^N}|\xi|^{2s}|\widehat u(\xi)|^2\,d\xi
       =\frac{C(N,s)}2\Js(u,u).
 \end{aligned}
\end{equation}
Thus $b(u)=2q(u)$ in the notation of the preceding subsections. Write
\[
 a=\frac{N(p-2)}2,\qquad d=a-2s>0,\qquad
\beta=p-\frac as\geq0.
\]
Here $\beta=0$ exactly when $p=2_s^*$.
For the mass-preserving dilation $u_t(x)=t^{N/2}u(tx)$,
\begin{equation}\label{am-fiber}
 \begin{split}
 A(u_t)&=t^{2s}\bigl(A(u)+2q(u)\ln t\bigr),\\
 I(u_t)&=t^{2s}\left(\frac12 A(u)+q(u)\ln t\right)
                       -\frac{t^a}{p}P(u),\\
 Q_I(u)&:=\left.\frac{d}{dt}I(u_t)\right|_{t=1}
                 =sA(u)+q(u)-\frac apP(u).
 \end{split}
\end{equation}
The logarithmic term is retained throughout the construction.

\subsubsection*{The radial Hilbert space and high-codimension barriers}
Let $X_{\rm rad}=H^{s+\log}_{\rm rad}(\R^N)$ over the real numbers and
\[
 S_c^{\rm rad}=\{u\in X_{\rm rad}:\|u\|_2=c\}.
\]
Choose $\omega>1/(es)$ sufficiently large that
\begin{equation}\label{am-norm}
 \|u\|_\omega^2:=A(u)+\omega\|u\|_2^2
\end{equation}
is an equivalent Hilbert norm on $\HH$. Here is a direct verification. The estimates \eqref{infinitesimal} and \eqref{tail} give constants $D_s,D_s'>0$ such that
\[
 \frac12K(u)-D_s\|u\|_2^2\leq A(u)
          \leq2K(u)+D_s'\|u\|_2^2.
\]
Taking $\omega>D_s+1$ proves the equivalence. The lower spectral bound is also immediate from
\[
 \min_{r\geq0}r^{2s}\ln r^2=-\frac1{es}.
\]
In particular, on $S_c^{\rm rad}$,
\begin{equation}\label{am-normenergy}
 I(u)=\frac12\|u\|_\omega^2-\frac\omega2c^2-\frac1p\|u\|_p^p.
\end{equation}

\begin{lemma}[Radial compactness at the Sobolev endpoint]\label{green-radial-endpoint}
For $N\geq2$ and $0<s<1$, the embedding
\begin{equation}\label{am-radialcompact}
 X_{\rm rad}\hookrightarrow L^r(\R^N)
 \quad\text{is compact for every }2<r\leq2_s^*.
\end{equation}
\end{lemma}
\begin{proof}
Let $w_n\rightharpoonup0$ in $X_{\rm rad}$. The inclusion $X\subset H^s$ and radial compactness give $w_n\to0$ in $L^{r_0}$ for any fixed $2<r_0<2_s^*$. To recall why radial compactness holds in that open range, local compactness excludes concentration in balls with bounded centers. If a ball with center tending to infinity retained a positive amount of $L^2$ mass, rotations in dimension $N\geq2$ would produce arbitrarily many disjoint balls with that same mass, contradicting the uniform $L^2$ bound. The fractional vanishing criterion then applies.

It remains to prove convergence at $r=2_s^*$. Choose a smooth radial cutoff $\chi$ with $0\leq\chi\leq1$, $\chi=1$ on $B_1$, and $\chi=0$ outside $B_2$. For $R>e$ define
\[
 \widehat{w_n^{\rm low}}(\xi)=\chi(\xi/R)\widehat w_n(\xi),
 \qquad w_n^{\rm high}=w_n-w_n^{\rm low}.
\]
The Fourier representation and boundedness in $X$ imply
\[
 \int_{|\xi|>R}|\xi|^{2s}|\widehat w_n|^2\,d\xi
 \leq\frac{1}{2\ln R}
       \int_{|\xi|>R}\ell(\xi)|\widehat w_n|^2\,d\xi
 \leq\frac{C}{\ln R}.
\]
For the last inequality, the negative part of $\ell$ is bounded by $1/(es)$, while $A(w_n)$ and $\|w_n\|_2$ are bounded. The fractional Sobolev inequality therefore gives
\[
 \|w_n^{\rm high}\|_{2_s^*}\leq C(\ln R)^{-1/2}.
\]
For fixed $R$, the low-frequency multiplier is convolution with a Schwartz kernel, so $\|w_n^{\rm low}\|_{r_0}\to0$. Fourier inversion also gives $\|w_n^{\rm low}\|_\infty\leq C_R\|w_n\|_2\leq C_R$. Interpolating these two bounds yields $\|w_n^{\rm low}\|_{2_s^*}\to0$. First let $n\to\infty$ and then $R\to\infty$ to obtain endpoint convergence. Interpolation gives the entire asserted range. No radial monotonicity or rearrangement is used.
\end{proof}

Choose nested finite-dimensional subspaces $V_k\subset X_{\rm rad}$ with $\dim V_k=k$ and dense union, and set $Z_k=V_{k-1}^{\perp_\omega}$. Define
\[
 \eta_k=\sup\{\|u\|_p:u\in Z_k,\ \|u\|_\omega=1\}.
\]
Then $\eta_k\to0$. Indeed, a bounded sequence of unit vectors in $Z_k$ converges weakly to zero, since it is eventually orthogonal to every fixed $V_j$; compactness in \eqref{am-radialcompact} gives convergence to zero in $L^p$.

For large $k$ put
\begin{equation}\label{am-radius}
 R_k=\left(\frac{p}{4\eta_k^p}\right)^{1/(p-2)},\qquad
 B_k=\{u\in S_c^{\rm rad}\cap Z_k:\|u\|_\omega=R_k\}.
\end{equation}
We may restrict to indices for which $R_k>2\sqrt\omega c$. By \eqref{am-normenergy},
\begin{equation}\label{am-barrier}
 \inf_{B_k}I\geq\frac{R_k^2}{4}-\frac\omega2c^2
                         \longrightarrow+\infty.
\end{equation}

\subsubsection*{The odd linking class}
The set $S_c^{\rm rad}\cap V_k$ is a compact sphere of dimension $k-1$. Formula \eqref{am-fiber} implies, uniformly on this sphere,
\[
 I(u_t)\to0\quad(t\downarrow0),\qquad
 I(u_t)\to-\infty\quad(t\to\infty),
\]
and
\[
 \|u_t\|_\omega\to\sqrt\omega c\quad(t\downarrow0),\qquad
 \|u_t\|_\omega\to\infty\quad(t\to\infty).
\]
For the uniform assertions, $P(u)$ and $q(u)$ have positive minima on the finite-dimensional sphere, while $A(u),P(u),q(u)$ are bounded there. In particular, $I(u_t)<0$ uniformly when $t$ is sufficiently small, because $q(u)\ln t$ is negative. Choose $0<t_k^-<t_k^+$ such that the endpoint maps
\[
 g_k^-(u)=u_{t_k^-},\qquad g_k^+(u)=u_{t_k^+}
\]
satisfy
\begin{equation}\label{am-endpoints}
 \|g_k^-(u)\|_\omega<R_k/2,\quad
 \|g_k^+(u)\|_\omega>2R_k,\quad
 I(g_k^\pm(u))\leq0.
\end{equation}
Let $\Gamma_k$ consist of continuous maps
\[
 g:[0,1]\times(S_c^{\rm rad}\cap V_k)\longrightarrow S_c^{\rm rad}
\]
that are odd in $u$ and equal $g_k^-$ and $g_k^+$ at the two endpoints. This class is nonempty: interpolate continuously between $\ln t_k^-$ and $\ln t_k^+$ in the dilation parameter.

Every $g\in\Gamma_k$ meets $B_k$ by the sphere-linking lemma of \cite[Lemma 2.3]{BartschValeriola}. For completeness, let $\Pi_{k-1}$ be the $\omega$-orthogonal projection onto $V_{k-1}$ and consider
\[
 (\tau,u)\longmapsto
 H_g(\tau,u)=\bigl(\Pi_{k-1}g(\tau,u),\|g(\tau,u)\|_\omega\bigr).
\]
The group $\{1,-1\}$ acts antipodally on the first coordinate and trivially on the second. The endpoint images lie on opposite sides of height $R_k$. The crossing property of the $\mathbb Z_2$ cohomological index implies that the image at height $R_k$ has index at least $k$, the index of $S_c^{\rm rad}\cap V_k$. If its first coordinate were never zero, projection and normalization would give an equivariant map into the sphere of $V_{k-1}$, whose index is $k-1$, a contradiction. This is the odd linking lemma; its use requires no sign assumption on the operator kernel.

Define the finite minimax levels
\begin{equation}\label{am-levels}
 h_k(c)=\inf_{g\in\Gamma_k}\max_{(\tau,u)}I(g(\tau,u)).
\end{equation}
The linking property and \eqref{am-barrier} yield
\begin{equation}\label{am-leveldiverge}
 h_k(c)\geq\inf_{B_k}I\longrightarrow+\infty.
\end{equation}

\subsubsection*{A Palais--Smale sequence with the scaling identity}
For each fixed sufficiently large $k$, write $h=h_k(c)>0$. Introduce the stretched functional
\begin{equation}\label{am-stretched}
 \widetilde I(\sigma,u)=I(u_{e^\sigma})
 =e^{2s\sigma}\left(\frac12A(u)+\sigma q(u)\right)
              -\frac{e^{a\sigma}}pP(u)
 \quad\text{on }\R\times S_c^{\rm rad}.
\end{equation}
This explicit scalar formula shows that $\widetilde I$ is of class $C^1$. The involution is $(\sigma,u)\mapsto(\sigma,-u)$. Lift $\Gamma_k$ to equivariant maps into $\R\times S_c^{\rm rad}$ with endpoint values $(0,g_k^-)$ and $(0,g_k^+)$. Projection by $(\sigma,u)\mapsto u_{e^\sigma}$ sends each lifted map into $\Gamma_k$; conversely, $g\mapsto(0,g)$ is a lift. The two minimax levels are therefore equal.

\begin{lemma}[Quantitative stretched minimax selection]\label{green-minimax}
For every sufficiently large fixed $k$, there exist $(\sigma_n,v_n)\in\R\times S_c^{\rm rad}$ such that
\[
 \widetilde I(\sigma_n,v_n)\to h_k(c),\qquad
 \bigl\|D(\widetilde I|_{\R\times S_c^{\rm rad}})(\sigma_n,v_n)\bigr\|\to0,
 \qquad \sigma_n\to0.
\]
\end{lemma}
\begin{proof}
Equip $\R\times S_c^{\rm rad}$ with the product metric induced by $\R\times X_{\rm rad}$. Choose $\epsilon_n\downarrow0$ with $h_k(c)-2\epsilon_n>0$ and $g_n\in\Gamma_k$ whose maximal energy is at most $h_k(c)+\epsilon_n$. Lift $g_n$ as $(0,g_n)$ and denote its compact image by $K_n$. The quantitative equivariant minimax deformation used in \cite[Lemmas 2.4--2.5]{BartschValeriola} selects a point at distance at most $2\sqrt{\epsilon_n}$ from $K_n$, with energy in $[h_k(c)-2\epsilon_n,h_k(c)+2\epsilon_n]$ and constrained derivative norm at most $C\sqrt{\epsilon_n}$, for an absolute constant $C$. To justify the selection, if no such point exists, an odd pseudo-gradient supported in this neighborhood and energy strip decreases energy at a rate proportional to $\sqrt{\epsilon_n}$. Running it for a time proportional to $\sqrt{\epsilon_n}$ lowers the maximum below $h_k(c)$, contradicting its definition. The deformation fixes the endpoints because their energies are nonpositive. Since every point of $K_n$ has first coordinate zero, the selected point satisfies $|\sigma_n|\leq2\sqrt{\epsilon_n}$. This gives all three conclusions.
\end{proof}

Set $u_n=(v_n)_{e^{\sigma_n}}$. Differentiation of \eqref{am-stretched} in $\sigma$ gives $Q_I(u_n)\to0$. Differentiation in $v$ and the uniformly bounded dilation operators for $\sigma_n\to0$ give
\begin{equation}\label{am-PS}
 \|u_n\|_2=c,\qquad I(u_n)\to h,\qquad Q_I(u_n)\to0,
 \qquad D(I|_{S_c^{\rm rad}})(u_n)\to0.
\end{equation}
This uses the chain rule for the scalar expression \eqref{am-stretched}, not differentiability of the dilation orbit in the energy norm.

\subsubsection*{Boundedness for every prescribed mass}
Two exact identities supply the estimates specific to the logarithmic operator. First,
\begin{equation}\label{am-quadraticidentity}
 I(u)-\frac1aQ_I(u)=\frac{d}{2a}A(u)-\frac1a q(u).
\end{equation}
Second, Jensen's inequality with the probability density $|\widehat u|^2/c^2$ gives
\begin{equation}\label{am-jensen}
 A(u)\geq\frac{q(u)}s\ln\frac{q(u)}{c^2}
       \qquad(u\in S_c^{\rm rad}).
\end{equation}
Indeed, apply the convexity of $y\ln y$ to $y=|\xi|^{2s}$. From \eqref{am-quadraticidentity}--\eqref{am-jensen},
\begin{equation}\label{am-logcoercivity}
 I(u)-\frac1aQ_I(u)
 \geq\frac{q(u)}a
       \left(\frac{d}{2s}\ln\frac{q(u)}{c^2}-1\right).
\end{equation}
The right-hand side tends to $+\infty$ as $q(u)\to\infty$. Applying this to \eqref{am-PS} proves that $q(u_n)$ is bounded. The exact identity
\[
 A(u_n)=\frac2d\bigl(aI(u_n)-Q_I(u_n)+q(u_n)\bigr)
\]
then proves that $A(u_n)$, and hence $\|u_n\|_\omega$, is bounded. No smallness assumption on $c$ enters this step.

\subsubsection*{The multiplier at high energy}
There are bounded multipliers $\lambda_n$ such that
\[
 \mathscr L u_n+\lambda_n u_n-|u_n|^{p-2}u_n\longrightarrow0
                  \quad\text{in }X_{\rm rad}^*.
\]
They may be chosen as $\lambda_n c^2=P(u_n)-A(u_n)$. Using \eqref{am-fiber} and \eqref{am-quadraticidentity}, this becomes the exact identity
\begin{equation}\label{am-multiplier}
 \lambda_n c^2
 =\frac{2s\beta}{d}I(u_n)
  +\frac{p-2}{d}\bigl(q(u_n)-Q_I(u_n)\bigr).
\end{equation}
After taking a subsequence, $\lambda_n\to\lambda$. If $p<2_s^*$, then $\beta>0$ and
\begin{equation}\label{am-gap}
 \lambda\geq\frac{2s\beta}{dc^2}h.
\end{equation}
In this subcritical case, choose $k$ large enough that
\begin{equation}\label{am-threshold}
 h_k(c)>\frac{dc^2}{2e s^2\beta}.
\end{equation}
For fixed $c>0$, this is possible by \eqref{am-leveldiverge} and gives $\lambda>1/(es)$.

At $p=2_s^*$ the preceding division by $\beta$ is unavailable. Instead, $a=sp$ and $d=s(p-2)$, so \eqref{am-multiplier} reduces to
\begin{equation}\label{green-critical-multiplier}
 \lambda_n c^2=\frac1s\bigl(q(u_n)-Q_I(u_n)\bigr).
\end{equation}
Let $C_S>0$ satisfy the fractional Sobolev estimate
\[
 P(v)\leq C_S q(v)^{p/2},\qquad p=2_s^*,
\]
and put $C_*=(p-2)C_S/(2p)$. Boundedness allows extraction with $q(u_n)\to q_\infty$. Since $Q_I(u_n)\to0$, \eqref{green-critical-multiplier} gives $\lambda c^2=q_\infty/s$. The exact identity
\[
 I(u_n)=\frac{p-2}{2p}P(u_n)-\frac12\lambda_n c^2
\]
therefore yields
\[
 h\leq C_*q_\infty^{p/2}-\frac{q_\infty}{2s}
 \leq C_*q_\infty^{p/2}.
\]
In particular,
\begin{equation}\label{green-critical-gap}
 \lambda\geq\frac1{sc^2}\left(\frac{h}{C_*}\right)^{2/p}.
\end{equation}
Choose $k$ so large that $h_k(c)>C_*(c^2/e)^{p/2}$. Then $\lambda>1/(es)$ also at the endpoint. Thus $\mathscr L+\lambda$ is coercive in both cases, and the following compactness argument applies throughout the enlarged exponent range.

\subsubsection*{Strong convergence and exact mass}
By boundedness and radial compactness, after extraction,
\[
 u_n\rightharpoonup u\quad\text{in }X_{\rm rad},\qquad
 u_n\to u\quad\text{in }L^p(\R^N).
\]
The power nonlinearity converges strongly in $L^{p/(p-1)}$, so the limit solves
\[
 \mathscr L u+\lambda u=|u|^{p-2}u
 \quad\text{against radial test functions}.
\]
Subtract the limit equation from the residual equation and test with $u_n-u$. We obtain
\[
 A(u_n-u)+\lambda\|u_n-u\|_2^2\longrightarrow0.
\]
For $\lambda>1/(es)$, this quadratic form controls the full energy norm. More explicitly, since $\ell(\xi)=|\xi|^{2s}\ln|\xi|^2\geq-1/(es)$,
\begin{equation}\label{am-coercivegap}
 A(v)+\lambda\|v\|_2^2\geq\kappa_{\lambda,\omega}\|v\|_\omega^2,
 \qquad
 \kappa_{\lambda,\omega}
 =\min\left\{1,\frac{\lambda-1/(es)}{\omega-1/(es)}\right\}>0.
\end{equation}
Thus $u_n\to u$ strongly in $H^{s+\log}$, and in particular
\[
 \|u\|_2=c,\qquad I(u)=h_k(c)>0,\qquad Q_I(u)=0.
\]
The limit cannot be zero because its norm is exactly $c$.

Finally, radial criticality implies criticality in the full space. For any test function, average it over the orthogonal group. At the radial function $u$, the derivative of $I+\lambda\|\cdot\|_2^2/2$ takes the same value on a test function and its radial average, by rotational invariance. Since it vanishes on radial tests, it vanishes on every test. Hence $u$ is a weak solution of the original stationary equation \eqref{equ3}.

Applying this argument for all sufficiently large $k$ proves Theorem~\ref{theorem-all-masses}. The solution energies tend to infinity; passing to a subsequence gives distinct pairs $\{u_{c,j},-u_{c,j}\}$. For $p<2_s^*$, formula \eqref{am-multiplier} implies $\lambda_{c,j}\to\infty$. At $p=2_s^*$, the endpoint estimate \eqref{green-critical-gap} gives the same conclusion. Their dilation curvature is strictly negative:
\[
 \left.\frac{d}{dt}Q_I((u_{c,j})_t)\right|_{t=1}
   =-d\,q(u_{c,j})-2sa\,I(u_{c,j})<0.
\]
They are therefore high-energy critical points on the physical upper Poho\u{z}aev branch. Their existence is compatible with a negative infimum over that branch at large mass, because they are not asserted to minimize that infimum.

\begin{remark}
The argument requires $N\geq2$ so that radial compactness applies. It does not establish the corresponding all-mass statement in dimension one. It also makes no positivity or ground-state assertion. In particular, no rearrangement inequality for the signed logarithmic energy is used.
\end{remark}

\appendix
\section*{Appendix}
\setcounter{equation}{0}
\renewcommand{\theequation}{A.\arabic{equation}}

\renewcommand{\thetheorem}{A.\arabic{theorem}}
\setcounter{theorem}{0}

\subsection*{Sharp interpolation and its optimizer}
In this part of the appendix, we assume $2\leq p\leq2_s^*$ when $N>2s$, and $2\leq p<\infty$ when $N\leq2s$. Thus the finite Sobolev endpoint is included. Define
\[
 T(u)=\|u\|_{\HH}^2=\|u\|_2^2+K(u),\qquad
 r_0=\frac{N(p-2)}{4s},\qquad
 \beta=p-\frac{N(p-2)}{2s}.
\]
The quotient associated with \eqref{equ13} is
\begin{equation}\label{quotient}
 J(u)=\frac{T(u)^{r_0}\|u\|_2^\beta}{\|u\|_p^p},\qquad u\ne0.
\end{equation}
For $p>2$, $r_0>0$ and $\beta\geq0$; equality $\beta=0$ occurs precisely at the finite endpoint $p=2_s^*$. Amplitude multiplication leaves this quotient unchanged; spatial dilation generally does not.

Indeed, the exact scaling of the positive logarithmic form is
\begin{equation}\label{Kscale}
 K(u_t)=Ct^{2s}\iint_{\R^{2N}}
 \frac{|u(X)-u(Y)|^2}{|X-Y|^{N+2s}}
       \bigl(\ln t-\ln|X-Y|\bigr)_+\,dX\,dY.
\end{equation}
The moving positive-part cutoff prevents the classical two-norm normalization. Translation preserves the norms, but cannot correct this dilation defect.

\begin{proposition}[Attainment of the interpolation quotient]\label{proposition1}
For $N\geq1$, $0<s<1$, and $2\leq p\leq2_s^*$ if $N>2s$ (or $2\leq p<\infty$ if $N\leq2s$), the quotient \eqref{quotient} attains its infimum over $\HH\setminus\{0\}$. A minimizer can be chosen nonnegative, radial, and nonincreasing. If $p=2$, $J(u)=1$ for every nonzero $u$, and $\mathbf K_{\rm opt}=1$. For $p>2$, $\mathbf K_{\rm opt}=(\inf J)^{-1}$.
\end{proposition}
\begin{proof}
We prove the nontrivial case $p>2$. The kernel of $K$ is nonnegative, radial, and nonincreasing.  Taking the absolute value and then the Schwarz rearrangement does not increase $K$ and preserves both the $L^2$-norm and the $L^p$-norm. The rearrangement inequality follows by truncating the kernel to an integrable radial decreasing kernel, expanding the quadratic difference, applying the Riesz inequality, and passing to the limit monotonically. Hence a minimizing sequence can be chosen nonnegative, radial, and nonincreasing. Using amplitude normalization alone, we may assume that
\[
 \|u_n\|_2=1.
\]
The continuous inclusion $\HH\subset H^s$ and the fractional Gagliardo--Nirenberg inequality (the fractional Sobolev inequality when $\beta=0$) show that
\[
 0<I:=\inf_{u\ne0}J(u)<\infty.
\]
We now prove that this normalized minimizing sequence is bounded in $\HH$.

\medskip\noindent\emph{Step 1: $T(u_n)$ is bounded.}
If a subsequence had $T_n:=T(u_n)\to\infty$, then \eqref{infinitesimal} and $K(u_n)=T_n-1$ would give, for every $\varepsilon>0$,
\[
 \limsup_{n\to\infty}\frac{S(u_n)}{T_n}\leq\varepsilon.
\]
Thus $S(u_n)/T_n\to0$. The same fractional interpolation inequality, including its Sobolev endpoint form when $\beta=0$, with $\|u_n\|_2=1$, implies
\[
 J(u_n)=\frac{T_n^{r_0}}{\|u_n\|_p^p}
 \geq\frac1{G_{N,s,p}}
       \left(\frac{T_n}{S(u_n)}\right)^{r_0}\longrightarrow\infty,
\]
contradicting $J(u_n)\to I$. This proves boundedness in $\HH$ without using a spatial dilation.

\medskip\noindent\emph{Step 2: global $L^p$ compactness and a nonzero limit.}
After passing to a subsequence, we have $u_n\rightharpoonup v$ in $\HH$. Proposition~\ref{proposition-embeddings}, including local compactness at the finite Sobolev endpoint, gives $u_n\to v$ in $L^p_{\rm loc}$ and almost everywhere. The limit remains nonnegative, radial, and nonincreasing. The radial monotonicity estimate gives
\[
 |u_n(x)|\leq\omega_N^{-1/2}|x|^{-N/2},
 \qquad \omega_N=|B_1|,
\]
because $\|u_n\|_2=1$. Therefore
\[
 \int_{|x|>R}|u_n|^p\,dx
 \leq\left(\omega_N^{-1/2}R^{-N/2}\right)^{p-2}
       \int_{|x|>R}|u_n|^2\,dx
 \leq\omega_N^{-(p-2)/2}R^{-N(p-2)/2}.
\]
The same tail bound holds for $v$. Combining it with local compactness yields $u_n\to v$ globally in $L^p$. Crucially, since $T_n\geq1$,
\begin{equation}\label{nonzero}
 \|u_n\|_p^p=\frac{T_n^{r_0}}{J(u_n)}\geq\frac1{J(u_n)}.
\end{equation}
Hence $\|v\|_p^p\geq1/I>0$ and $v\ne0$. Strong $L^p$ convergence by itself would not imply this nonvanishing; estimate \eqref{nonzero} is indispensable.

\medskip\noindent\emph{Step 3: attainment.}
Choose a further subsequence with $T_n\to T_*<\infty$. Weak lower semicontinuity gives $T(v)\leq T_*$ and $\|v\|_2\leq1$. Since $\newadd{\beta\geq0}$,
\[
 I\leq J(v)=\frac{T(v)^{r_0}\|v\|_2^\beta}{\|v\|_p^p}
 \leq\frac{T_*^{r_0}}{\|v\|_p^p}
 =\lim_{n\to\infty}J(u_n)=I.
\]
Thus $v$ is an optimizer. Since $r_0>0$, equality in this chain forces $T(v)=T_*$. Convergence of the Hilbert norms together with weak convergence gives strong convergence in $\HH$, and hence $\|v\|_2=1$. This argument also covers $\beta=0$ without inferring mass recovery from a factor whose exponent vanishes.
\end{proof}

\begin{remark}[The attained endpoint embedding constant]
When $N>2s$ and $p=2_s^*$, Proposition~\ref{proposition1} gives
\[
 \|u\|_{2_s^*}^2\leq C_{\log}\|u\|_{\HH}^2,
 \qquad
 C_{\log}=\left(\inf_{u\ne0}J(u)\right)^{-2/2_s^*},
\]
with equality for a nonzero nonnegative radial nonincreasing function. Here the norm uses the positive logarithmic form $K$; this conclusion does not assert rearrangement monotonicity for the full signed form $A$.
\end{remark}

\subsection*{Poho\u{z}aev identity and admissibility of dilation tests}
\begin{lemma}\label{lemma2}
Let $2<p<2_s^*$ and let $u\in\Sc$ be a constrained critical point of $\E$. If the curve $t\mapsto u_t$ is differentiable at $t=1$ in $\HH$, then
\[
 2s\bigl(\Jp(u,u)-\Jm(u,u)\bigr)
 +C(N,s)\bigl(1+sB(N,s)\bigr)\Js(u,u)
 =\frac{N(p-2)}{2p}\|u\|_p^p.
\]
The same identity holds for every translated weak profile of a bounded minimizing sequence from Theorem~\ref{theorem4}, without the differentiability assumption on its dilation orbit.
\end{lemma}
\begin{proof}
In the first case, the dilation curve lies in $\Sc$, so constrained criticality and the chain rule imply
\[
 0=\left.\frac{d}{dt}\E(u_t)\right|_{t=1}=\Q(u).
\]
Substituting the explicit fiber formula gives the displayed identity. In the second case, apply Lemma~\ref{profilevariation}: replacing only the weak profile by its dilation, normalizing the total mass, and projecting back to the upper branch yields $\Q(u)=0$ by scalar differentiation. In particular, the proof of Theorem~\ref{theorem4} requires no unproved Poho\u{z}aev identity for arbitrary weak critical points.
\end{proof}

\begin{center}
	\textbf{Statements and Declarations}
\end{center}
\textbf{Ethics approval and consent to participate.} Not applicable.\\[4pt]
\textbf{Funding.} Not applicable.\\[4pt]
\textbf{Availability of data and materials.} Not applicable.\\[4pt]
\textbf{Conflict of interest.} The authors declare that there is no conflict of interest.

\addtocontents{toc}{\protect\setcounter{tocdepth}{2}}

\end{document}